\documentclass[english,romanappendices]{IEEEtran}
\usepackage[T1]{fontenc}
\usepackage[latin9]{inputenc}
\usepackage{verbatim}
\usepackage{float}
\usepackage{amsthm}
\usepackage{amstext}
\usepackage{amssymb}
\usepackage{graphicx}

\makeatletter

\floatstyle{ruled}
\newfloat{algorithm}{tbp}{loa}
\providecommand{\algorithmname}{Algorithm}
\floatname{algorithm}{\protect\algorithmname}

\theoremstyle{plain}
\newtheorem{thm}{\protect\theoremname}
\theoremstyle{definition}
\newtheorem{defn}[thm]{\protect\definitionname}
\theoremstyle{definition}
\newtheorem{problem}[thm]{\protect\problemname}
\theoremstyle{plain}
\newtheorem{cor}[thm]{\protect\corollaryname}
\theoremstyle{plain}
\newtheorem{prop}[thm]{\protect\propositionname}
\theoremstyle{plain}
\newtheorem{lem}[thm]{\protect\lemmaname}
\theoremstyle{remark}
\newtheorem{rem}[thm]{\protect\remarkname}

\usepackage[caption=false,font=footnotesize]{subfig}
\usepackage{url}

\usepackage{amsmath}
\usepackage{amsfonts}

\newcommand{\norm}[1]{\left\Vert#1\right\Vert}

\newcommand{\optmin}{\mathrm{min.}}

\newcommand{\optst}{\;\mathrm{s.t.}}

\allowdisplaybreaks[3]

\newcommand{\allthanks}{\thanks{The authors are with the Department of Electrical and Systems Engineering, University of Pennsylvania, Philadelphia, PA 19104. {\tt\scriptsize {\{hanshuo,utopcu,pappasg\}@seas.upenn.edu}} This work was supported in part by the NSF (CNS-1239224) and TerraSwarm, one of six centers of STARnet, a Semiconductor Research Corporation program sponsored by MARCO and DARPA.}}

\makeatother

\usepackage{babel}
\providecommand{\corollaryname}{Corollary}
\providecommand{\definitionname}{Definition}
\providecommand{\lemmaname}{Lemma}
\providecommand{\problemname}{Problem}
\providecommand{\propositionname}{Proposition}
\providecommand{\remarkname}{Remark}
\providecommand{\theoremname}{Theorem}

\begin{document}

\title{\textbf{\Large{}Differentially Private Distributed Constrained Optimization}}

\author{Shuo Han, Ufuk Topcu, George J. Pappas\allthanks}
\maketitle
\begin{abstract}
Many resource allocation problems can be formulated as an optimization
problem whose constraints contain sensitive information about participating
users. This paper concerns solving this kind of optimization problem
in a distributed manner while protecting the privacy of user information.
Without privacy considerations, existing distributed algorithms normally
consist in a central entity computing and broadcasting certain public
coordination signals to participating users. However, the coordination
signals often depend on user information, so that an adversary who
has access to the coordination signals can potentially decode information
on individual users and put user privacy at risk. We present a distributed
optimization algorithm that preserves differential privacy, which
is a strong notion that guarantees user privacy regardless of any
auxiliary information an adversary may have. The algorithm achieves
privacy by perturbing the public signals with additive noise, whose
magnitude is determined by the sensitivity of the projection operation
onto user-specified constraints. By viewing the differentially private
algorithm as an implementation of stochastic gradient descent, we
are able to derive a bound for the suboptimality of the algorithm.
We illustrate the implementation of our algorithm via a case study
of electric vehicle charging. Specifically, we derive the sensitivity
and present numerical simulations for the algorithm. Through numerical
simulations, we are able to investigate various aspects of the algorithm
when being used in practice, including the choice of step size, number
of iterations, and the trade-off between privacy level and suboptimality.

\end{abstract}

\section{Introduction}

Electric vehicles (EVs), including pure electric and hybrid plug-in
vehicles, are believed to be an important component of future power
systems~\cite{ipakchi2009grid}. Studies predict that the market
share of EVs in the United States can reach approximately 25\% by
year 2020~\cite{hadley2009potential}. By that time, EVs will become
a significant load on the power grid~\cite{taylor2010evaluations,clement2010impact},
which can lead to undesirable effects such as voltage deviations if
charging of the vehicles are uncoordinated. 

The key to reducing the impact of EVs on the power grid is to coordinate
their charging schedules, which is often cast as a constrained optimization
problem with the objective of minimizing the peak load, power loss,
or load variance~\cite{sortomme2011coordinated,deilami2011real}.
Due to the large number of vehicles, computing an optimal schedule
for all vehicles can be very time consuming if the computation is
carried out on a centralized server that collects demand information
from users. Instead, it is more desirable that the computation is
distributed to individual users. Among others, Ma et al.~\cite{ma2013decentralized}
proposed a distributed charging strategy based on the notion of valley-filling
charging profiles, which is guaranteed to be optimal when all vehicles
have identical (i.e., homogeneous) demand. Gan et al.~\cite{gan2013optimal}
proposed a more general algorithm that is optimal for nonhomogeneous
demand and allows asynchronous communication.

In order to solve the constrained optimization problem of scheduling
in a distributed manner, the server is required to publish certain
public information that is computed based on the tentative demand
collected from participating users. Charging demand often contains
private information of the users. As a simple example, zero demand
from a charging station attached to a single home unit is a good indication
that the home owner is away from home. Note that the public coordination
signal is received by everyone including potential adversaries whose
goal is to decode private user information from the public signal,
so that it is desirable to develop solutions for protecting user privacy.

It has been long recognized that \emph{ad hoc} solutions such as anonymization
of user data are inadequate to guarantee privacy due to the presence
of public side information. A famous case is the reidentification
of certain users from an anonymized dataset published by Netflix,
which is an American provider of on-demand Internet streaming media.
The dataset was provided for hosting an open competition called the
Netflix Prize for find the best algorithm to predict user ratings
on films. It has been reported that certain Netflix subscribers can
be identified from the anonymized Netflix prize dataset through auxiliary
information from the Internet Movie Database (IMDb)~\cite{narayanan2007break}.
As such, providing rigorous solutions to preserving privacy has become
an active area of research. In the field of systems and control, recent
work on privacy includes, among others, filtering of streaming data~\cite{leny2014differentially},
smart metering~\cite{sankar2013smart}, traffic monitoring~\cite{canepa2013framework},
and privacy in stochastic control~\cite{venkitasubramaniam2013privacy}.

Recently, the notion of \emph{differential privacy} proposed by Dwork
and her collaborators has received attention due to its mathematically
rigorous formulation~\cite{dwork2006calibrating}. The original setting
assumes that the sensitive user information is held by a trustworthy
party (often called \emph{curator} in related literature), and the
curator needs to answer external queries (about the sensitive user
information) that potentially come from an adversary who is interested
in learning information belonging to some user. For example, in EV
charging, the curator is the central server that aggregates user information,
and the queries correspond to public coordination signals. Informally,
preserving differential privacy requires that the curator must ensure
that the results of the queries remain approximately unchanged if
data belonging to any single user are modified. In other words, the
adversary should know little about any single user's information from
the results of queries. A recent survey on differential privacy can
be found in~\cite{dwork2008differential}; there is also a recent
textbook on this topic written by Dwork and Roth~\cite{dwork2013algorithmic}.

\paragraph*{Contributions}

Motivated by the privacy concerns in EV charging and recent advances
in differential privacy, in this paper, we investigate the problem
of \emph{preserving differential privacy in distributed constrained
optimization}. We present a differentially private distributed algorithm
for solving a class of constrained optimization problems, whose privacy
guarantee is proved using the adaptive composition theorem. We show
that the private optimization algorithm can be viewed as an implementation
of stochastic gradient descent~\cite{robbins1951stochastic}. Based
on previous results on stochastic gradient descent~\cite{shamir2013stochastic},
we are able to derive a bound for the suboptimality of our algorithm
and reveal the trade-off between privacy and performance of the algorithm. 

We illustrate the applicability of this general framework of differentially
private distributed constrained optimization in the context of EV
charging. To this end, we begin by computing the \emph{sensitivity}
of the public signal with respect to changes in private information.
Specifically, this requires analyzing the sensitivity of the projection
operation onto the user-specified constraints. Although such sensitivity
can be difficult to compute for a general problem, using tools in
optimization theory, we are able to derive an explicit expression
of the sensitivity for the EV charging example. %
Through numerical simulations, we show that our algorithm is able
to provide strong privacy guarantees with little loss in performance
when the number of participating users (i.e., vehicles) is large.

\paragraph*{Related work}

There is a large body of research work on incorporating differential
privacy into resource allocation problems. A part of the work deals
with indivisible resources (or equivalently, games with discrete actions),
including the work by, among others, Kearns et al.~\cite{kearns2014mechanism},
Rogers and Roth~\cite{rogers2014asymptotically}, and Hsu et al.~\cite{hsu2013private}.
Our paper focuses on the case of \emph{divisible} resources and where
private information is contained in the \emph{constraints} of the
allocation problem. 

In the work of differentially private resource allocation, it is a
common theme that the coordination signals are randomly perturbed
to avoid revealing private information of the users, such as in the
work by Huang et al.~\cite{huang2014differentially} and Hsu et al.~\cite{hsu2014privately}.
Huang et al.~\cite{huang2014differentially} study the problem of
differentially private convex optimization \emph{in the absence of
constraints}. In their formulation, the private user information is
encoded in the individual cost functions and can be interpreted as
user preferences. Apart from incorporating constraints, a major difference
of our setting compared to Huang et al.~\cite{huang2014differentially}
is the way that privacy is incorporated into the optimization problem.
In Huang et al.~\cite{huang2014differentially}, they assume that
the public coordination signals \emph{do not} change when the individual
cost function changes. However, this assumption fails to hold, e.g.,
in the case of EV charging, where the coordination signal is computed
by aggregating the response from users and hence is sensitive to changes
in user information. Instead, we treat the public coordination signals
as the quantity that needs to be perturbed (i.e., as a \emph{query},
in the nomenclature of differential privacy) in order to prevent privacy
breach caused by broadcasting the signals. 

The recent work by Hsu et al.~\cite{hsu2014privately} on privately
solving linear programs is closely related to our work, since their
setting also assumes that the private information is contained in
the (affine) constraints. Our work can be viewed as a generalization
of their setting by extending the form of objective functions and
constraints. In particular, the objective function can be any convex
function that depends on the aggregate allocation and has Lipschitz
continuous gradients. The constraints can be any convex and separable
constraints; for illustration, we show how to implement the algorithm
for a particular set of affine constraints motivated by EV charging.

\paragraph*{Paper organization}

The paper is organized as follows. Section~\ref{sec:Background-material}
introduces the necessary background on (non-private) distributed optimization
and, in particular, projected gradient descent. Section~\ref{sec:Problem-formulation-and}
reviews the results in differential privacy and gives a formal problem
statement of differentially private distributed constrained optimization.
Section~\ref{sec:overview_of_results} gives an overview of the main
results of the paper. Section~\ref{sec:A-differentially-private}
describes a differentially private distributed algorithm that solves
a general class of constrained optimization problems. We also study
the trade-off between privacy and performance by analyzing the suboptimality
of the differentially private algorithm. 

In Section~\ref{sec:Sensitivity-computation}, we illustrate the
implementation of our algorithm via a case study of EV charging. In
particular, we compute the sensitivity of the projection operation
onto user-specified constraints, which is required for implementing
our private algorithm. Section~\ref{sec:Numerical-simulations} presents
numerical simulations on various aspects of the algorithm when being
used in practice, including choice of step size, number of iterations,
and the trade-off between privacy level and performance.

\section{Background: Distributed constrained optimization\label{sec:Background-material}}

\subsection{Notation}

Denote the $\ell_{p}$-norm of any~$x\in\mathbb{R}^{n}$ by~$\norm{x}_{p}$.
The subscript~$p$ is dropped in the case of~$\ell_{2}$-norm. For
any nonempty convex set~$\mathcal{C}\subset\mathbb{R}^{n}$ and $x\in\mathbb{R}^{n}$,
denote by~$\Pi_{\mathcal{C}}(x)$ the projection operator that projects
$x$ onto $\mathcal{C}$ in $\ell_{2}$-norm. Namely, $\Pi_{\mathcal{C}}(x)$
is the solution of the following constrained least-squares problem
\begin{equation}
\underset{\hat{x}}{\optmin}\quad\Vert\hat{x}-x\Vert^{2}\qquad\optst\quad\hat{x}\in\mathcal{C}.\label{eq:ls_projection}
\end{equation}
It can be shown that problem~(\ref{eq:ls_projection}) is always
feasible and has a unique solution so that~$\Pi_{\mathcal{C}}$ is
well-defined. For any function~$f$ (not necessarily convex), denote
by~$\partial f(x)$ the set of subgradients of~$f$ at $x$:
\[
\partial f(x):=\{g\colon f(y)\geq f(x)+g^{T}(y-x)\mbox{ for all }y\}.
\]
When $f$ is convex and differentiable at $x$, the set~$\partial f(x)$
becomes a singleton set whose only element is the gradient $\nabla f(x)$.
For any function $f$, denote its range by $\mathrm{range}(f)$. For
any differentiable function $f$ that depends on multiple variables
including $x$, denote by $\partial_{x}f$ the partial derivative
of $f$ with respect to $x$. For any~$\lambda>0$, denote by~$\mathrm{Lap}(\lambda)$
the zero-mean Laplace probability distribution such that the probability
density function of a random variable $X$ obeying the distribution
$\mathrm{Lap}(\lambda)$ is~$p_{X}(x)=\frac{1}{2\lambda}\exp(-|x|/\lambda)$.
The vector consisting all ones is written as~$\mathbf{1}$. The symbol~$\preceq$
is used to represent element-wise inequality: for any~$x,y\in\mathbb{R}^{n}$,
we have~$x\preceq y$ if and only if~$x_{i}\leq y_{i}$ for all~$1\leq i\leq n$.
For any positive integer $n$, we denote by $[n]$ the set $\{1,2,\dots,n\}$.

\subsection{Distributed constrained optimization}

Before discussing privacy issues, we first introduce the necessary
background on distributed constrained optimization. We consider a
constrained optimization problem over $n$ variables $r_{1},r_{2},\dots,r_{n}\in\mathbb{R}^{T}$
in the following form:
\begin{alignat}{2}
 & \underset{\{r_{i}\}_{i=1}^{n}}{\optmin}\quad &  & U\left({\textstyle \sum_{i=1}^{n}}r_{i}\right)\label{eq:dco_prob}\\
 & \underset{\phantom{\{r_{i}\}_{i=1}^{n}}}{\optst}\quad &  & r_{i}\in\mathcal{C}_{i},\quad i\in[n].\nonumber 
\end{alignat}
Throughout the paper, we assume that the objective function~$U\colon\mathbb{R}^{T}\to\mathbb{R}$
in problem~(\ref{eq:dco_prob}) is differentiable and convex, and
its gradient $\nabla U$ is $L$-Lipschitz in the $\ell_{2}$-norm,
i.e., there exists $L>0$ such that 
\[
\norm{\nabla U(x)-\nabla U(y)}\leq L\norm{x-y}\quad\text{for all }x,y.
\]
The set $\mathcal{C}_{i}$ is assumed to be convex for all $i\in[n]$.
For resource allocation problems, the variable $r_{i}$ and the constraint
set $\mathcal{C}_{i}$ are used to capture the allocation and constraints
on the allocation for user/agent~$i$.

\begin{algorithm}
\textbf{Input}: $U$, $\{\mathcal{C}_{i}\}_{i=1}^{n}$, $K$, and
step sizes $\{\alpha_{k}\}_{k=1}^{K}$.

\textbf{Output}: $\{r_{i}^{(K+1)}\}_{i=1}^{n}$.

Initialize $\{r_{i}^{(1)}\}_{i=1}^{n}$ arbitrarily. For $k=1,2,\dots,K$,
repeat:
\begin{enumerate}
\item Compute $p^{(k)}:=\nabla U\left({\textstyle \sum_{i=1}^{n}r_{i}^{(k)}}\right)$.
\item For $i\in[n]$, update~$r_{i}^{(k+1)}$ according to 
\begin{equation}
r_{i}^{(k+1)}:=\Pi_{\mathcal{C}_{i}}(r_{i}^{(k)}-\alpha_{k}p^{(k)}).\label{eq:ev_charging_prob}
\end{equation}
 %

\end{enumerate}
\caption{Distributed projected gradient descent (with a fixed number of iterations).\label{alg:Distributed-EV-charging}}
\end{algorithm}

The optimization problem~(\ref{eq:dco_prob}) can be solved iteratively
using projected gradient descent, which requires computation of the
gradient of~$U$ and its projection onto the feasible set at each
iteration. The computational complexity of projected gradient descent
is dominated by the projection operation and grows with~$n$. For
practical applications, the number~$n$ can be quite large, so that
it is desirable to distribute the projection operation to individual
users. A distributed version of the projected gradient descent method
is shown in Algorithm~\ref{alg:Distributed-EV-charging}. The algorithm
guarantees that the output converges to the optimal solution as~$K\to\infty$
with proper choice of step sizes~$\{\alpha_{k}\}_{k=1}^{K}$ (see~\cite{gan2013optimal}
for details on how to choose~$\alpha_{k}$).

\section{Problem formulation\label{sec:Problem-formulation-and}}

\subsection{Privacy in distributed constrained optimization}

In many applications, the specifications of $\mathcal{C}_{i}$ may
contain sensitive information that user $i$ wishes to keep undisclosed
from the public. %
In the framework of differential privacy, it is assumed that an adversary
can potentially collaborate with some users in the database in order
to learn about other user's information. Under this assumption, the
distributed projected descent algorithm (Algorithm~\ref{alg:Distributed-EV-charging})
can lead to possible loss of privacy of participating users for reasons
described below. It can be seen from Algorithm~\ref{alg:Distributed-EV-charging}
that $\mathcal{C}_{i}$ affects~$r_{i}^{(k)}$ through equation~(\ref{eq:ev_charging_prob})
and consequently also $p^{(k)}$. Since $p^{(k)}$ is broadcast publicly
to every charging station, with enough side information (such as collaborating
with some participating users), an adversary who is interested in
learning private information about some user~$i$ may be able to
infer information about $\mathcal{C}_{i}$ from the public signals~$\{p^{(k)}\}_{k=1}^{K}$.
We will later illustrate the privacy issues in the context of EV charging.

\subsection{Differential privacy}

Our goal is to modify the original distributed projected gradient
descent algorithm (Algorithm~\ref{alg:Distributed-EV-charging})
to preserve \emph{differential privacy}. Before giving a formal statement
of our problem, we first present some preliminaries on differential
privacy. Differential privacy considers a set (called \emph{database})
$D$ that contains private user information to be protected. For convenience,
we denote by $\mathcal{D}$ the universe of all possible databases
of interest. The information that we would like to obtain from a database~$D$
is given by $q(D)$ for some mapping~$q$ (called \emph{query}) that
acts on $D$. In differential privacy, preserving privacy is equivalent
to hiding changes in the database. Formally, changes in a database
can be defined by a symmetric binary relation between two databases
called \emph{adjacency} relation, which is denoted by~$\mathrm{Adj}(\cdot,\cdot)$;
two databases $D$ and $D'$ that satisfy $\mathrm{Adj}(D,D')$ are
called adjacent databases.
\begin{defn}
[Adjacent databases]Two databases $D=\{d_{i}\}_{i=1}^{n}$ and $D'=\{d_{i}'\}_{i=1}^{n}$
are said to be \emph{adjacent }if there exists $i\in[n]$ such that
$d_{j}=d_{j}'$ for all $j\neq i$. 
\end{defn}
A \emph{mechanism} that acts on a database is said to be differentially
private if it is able to ensure that two adjacent databases are nearly
indistinguishable from the output of the mechanism.
\begin{defn}
[Differential privacy~\cite{dwork2006calibrating}]\label{def:dp}Given
$\epsilon\geq0$, a mechanism~$M$ preserves $\epsilon$-differential
privacy if for all~$\mathcal{R}\subseteq\mathrm{range}(M)$ and all
adjacent databases~$D$ and~$D'$ in $\mathcal{D}$, it holds that
\begin{equation}
\mathbb{P}(M(D)\in\mathcal{R})\leq e^{\epsilon}\mathbb{P}(M(D')\in\mathcal{R}).\label{eq:dp_def}
\end{equation}

\end{defn}
The constant~$\epsilon$ indicates the level of privacy: smaller~$\epsilon$
implies higher level of privacy. The notion of differential privacy
promises that an adversary cannot tell from the output of $M$ with
high probability whether data corresponding to a single user in the
database have changed. It can be seen that any non-constant differentially
mechanism is necessarily \emph{randomized}, i.e., for a given database,
the output of such a mechanism obeys a certain probability distribution.
Finally, although it is not explicitly mentioned in Definition~\ref{def:dp},
a mechanism needs to be an approximation of the query of interest
in order to be useful. For this purpose, a mechanism is normally defined
in conjunction with some query of interest; a common notation is to
include the query~$q$ of interest in the subscript of the mechanism
as $M_{q}$.

\subsection{Problem formulation: Differentially private distributed constrained
optimization}

Recall that our goal of preserving privacy in distributed optimization
is to protect the user information in $\mathcal{C}_{i}$, even if
an adversary can collect all public signals $\{p^{(k)}\}_{k=1}^{K}$.
To mathematically formulate our goal under the framework of differential
privacy, we define the database~$D$ as the set $\{\mathcal{C}_{i}\}_{i=1}^{n}$
and the query as the~$K$-tuple consisting of all the gradients $p=(p^{(1)},p^{(2)},\dots,p^{(K)})$.
We assume that $\mathcal{C}_{1},\mathcal{C}_{2},\dots,\mathcal{C}_{n}$
belong a family of sets parameterized by $\alpha\in\mathbb{R}^{s}$.
Namely, there exists a parameterized set $\mathcal{C}$ such that
for all $i\in[n]$, we can write $\mathcal{C}_{i}=\mathcal{C}(\alpha_{i})$
for some $\alpha_{i}\in\mathbb{R}^{s}$. We also assume that there
exists a metric $\rho\colon\mathbb{R}^{s}\times\mathbb{R}^{s}\to\mathbb{R}_{+}$.
In this way, we can define the distance $\rho_{\mathcal{C}}(\mathcal{C}_{i},\mathcal{C}_{i}')$
between any $\mathcal{C}_{i}=\mathcal{C}(\alpha_{i})$ and $\mathcal{C}_{i}'=\mathcal{C}(\alpha_{i}')$
using the metric $\rho$ as 
\[
\rho_{\mathcal{C}}(\mathcal{C}_{i},\mathcal{C}_{i}'):=\rho(\alpha_{i},\alpha_{i}').
\]
For any given $\delta\mathcal{C}\in\mathbb{R}_{+}$, we define and
use throughout the paper the following adjacency relation between
any two databases $D$ and $D'$ in the context of distributed constrained
optimization. 
\begin{defn}
[Adjacency relation for constrained optimization]\label{def:adj_dco}For
any databases~$D=\{\mathcal{C}_{i}\}_{i=1}^{n}$ and~$D'=\{\mathcal{C}_{i}'\}_{i=1}^{n}$,
it holds that $\mathrm{Adj}(D,D')$ if and only if there exists $i\in[n]$
such that $\rho_{\mathcal{C}}(\mathcal{C}_{i},\mathcal{C}_{i}')\leq\delta\mathcal{C}$,
and $\mathcal{C}_{j}=\mathcal{C}_{j}'$ for all $j\neq i$.
\end{defn}
The constant $\delta\mathcal{C}$ is chosen based on the privacy requirement,
i.e., the kind of user activities that should be kept private. Using
the adjacency relation described in Definition~\ref{def:adj_dco},
we state in the following the problem of designing a differentially
private distributed algorithm for constrained optimization.
\begin{problem}
[Differentially private distributed constrained optimization]\label{prob:Find-a-randomized}Find
a randomized mechanism~$M_{p}$ that approximates the gradients~$p=(p^{(1)},p^{(2)},\dots,p^{(K)})$
(defined in Algorithm~\ref{alg:Distributed-EV-charging}) and preserves
$\epsilon$-differential privacy under the adjacency relation described
in Definition~\ref{def:adj_dco}. Namely, for any adjacent databases
$D$ and $D'$, and any $\mathcal{R}\subseteq\mathrm{range}(M_{p})$,
the mechanism $M_{p}$ should satisfy
\[
\mathbb{P}(M_{p}(D)\in\mathcal{R})\leq e^{\epsilon}\mathbb{P}(M_{p}(D')\in\mathcal{R}).
\]

\end{problem}

\subsection{Example application: EV charging\label{sub:Example-application:-EV}}

In EV charging, the goal is to charge $n$ vehicles over a horizon
of $T$ time steps with minimal influence on the power grid. For simplicity,
we assume that each vehicle belongs to one single user. For any $i\in[n]$,
the vector $r_{i}\in\mathbb{R}^{T}$ represents the charging rates
of vehicle~$i$ over time. In the following, we will denote by $r_{i}(t)$
the $t$-th component of $r_{i}$. Each vehicle needs to be charged
a given amount of electricity~$E_{i}>0$ by the end of the scheduling
horizon; in addition, for any $t\in[T]$, the charging rate $r_{i}(t)$
cannot exceed the maximum rate $\bar{r}_{i}(t)$ for some given constant
vector $\bar{r}_{i}\in\mathbb{R}^{T}$. Under these constraints on
$r_{i}$, the set $\mathcal{C}_{i}$ is described as follows: 
\begin{equation}
0\preceq r_{i}\preceq\bar{r}_{i},\qquad\mathbf{1}^{T}r_{i}=E_{i}.\label{eq:ev_charging_constr}
\end{equation}
The tuple $(\bar{r}_{i},E_{i})$ is called the \emph{charging specification}
of user $i$. Throughout the paper, we assume that $\bar{r}_{i}$
and $E_{i}$ satisfy 
\begin{equation}
\mathbf{1}^{T}\bar{r}_{i}\geq E_{i}\quad\text{for all }i\in[n],\label{eq:feasiblity}
\end{equation}
so that the constraints~(\ref{eq:ev_charging_constr}) are always
feasible. 

The objective function~$U$ in problem~(\ref{eq:dco_prob}) quantifies
the influence of a charging schedule~$\{r_{i}\}_{i=1}^{n}$ on the
power grid. We choose~$U$ as follows for the purpose of minimizing
load variance:
\begin{equation}
U\left({\textstyle \sum_{i=1}^{n}r_{i}}\right)=\frac{1}{2}\norm{d+{\textstyle \sum_{i=1}^{n}r_{i}}/m}^{2}.\label{eq:cost_quad}
\end{equation}
In~(\ref{eq:cost_quad}), $m$ is the number of households, which
is assumed proportional to the number of EVs, i.e., there exists $\gamma$
such that $n/m=\gamma$; then, the quantity ${\textstyle \sum_{i=1}^{n}r_{i}}/m$
becomes the aggregate EV load per household. The vector $d\in\mathbb{R}^{T}$
is the base load profile incurred by loads in the power grid other
than EVs, so that $U\left({\textstyle \sum_{i=1}^{n}r_{i}}\right)$
quantifies the variation of the total load including the base load
and EVs. It can be verified that $U$ is convex and differentiable,
and $\nabla U$ is Lipschitz continuous. 

The set $\mathcal{C}_{i}$ (defined by~$\bar{r}_{i}$ and $E_{i}$)
can be associated with personal activities of the owner of vehicle~$i$
in the following way. For example, $\bar{r}_{i}(t)=0$ may indicate
that the owner is temporarily away from the charging station (which
may be co-located with the owner's residence) so that the vehicle
is not ready to be charged. Similarly, $E_{i}=0$ may indicate that
the owner is not actively using the vehicle so that the vehicle does
not need to be charged. 

We now illustrate why publishing the exact gradient~$p^{(k)}$ can
potentially lead to a loss of privacy. The gradient $p^{(k)}$ can
be computed as $p^{(k)}=\frac{1}{m}(d+{\textstyle \sum_{i=1}^{n}r_{i}}/m)$.
If an adversary collaborates with all but one user $i$ so that the
adversary is able obtain $r_{j}^{(k)}$ for all $j\neq i$ in the
database. Then, the adversary can infer $r_{i}^{(k)}$ exactly from
$p^{(k)}$, even though user $i$ did not reveal his $r_{i}^{(k)}$
to the adversary. After obtaining $r_{i}^{(k)}$, the adversary can
obtain information on $\mathcal{C}_{i}$ by, for example, computing
$E_{i}=\mathbf{1}^{T}r_{i}^{(k)}$. 

The adjacency relation in the case of EV charging is defined as follows.
Notice that, in the case of EV charging, the parameter $\alpha_{i}$
that parameterizes the set $\mathcal{C}_{i}$ is given by $\alpha_{i}=(\bar{r}_{i},E_{i})$,
in which $(\bar{r}_{i},E_{i})$ is the charging specifications of
user $i$ as defined in~(\ref{eq:ev_charging_constr}). 
\begin{defn}
[Adjacency relation for EV charging]\label{def:adj_ev}For any databases
$D=\{\mathcal{C}_{i}(\bar{r}_{i},E_{i})\}_{i=1}^{n}$ and $D'=\{\mathcal{C}_{i}'(\bar{r}_{i}',E_{i}')\}_{i=1}^{n}$,
we have $\mathrm{Adj}(D,D')$ if and only if there exists $i\in[n]$
such that 
\begin{equation}
\norm{\bar{r}_{i}-\bar{r}_{i}'}_{1}\leq\delta r,\qquad|E_{i}-E_{i}'|\leq\delta E,\label{eq:adj_ev}
\end{equation}
and $\bar{r}_{j}=\bar{r}_{j}'$, $E_{j}=E_{j}'$ for all $j\neq i$. 
\end{defn}
In terms of choosing $\delta E$ and $\delta r$, one useful choice
for $\delta E$ is the maximum amount of energy an EV may need; this
choice of $\delta E$ can be used to hide the event corresponding
to whether an user needs to charge his vehicle.

\section{Overview of main results\label{sec:overview_of_results}}

\subsection{Results for general constrained optimization problems}

\begin{algorithm}
\textbf{Input}: $U$, $L$, $\{\mathcal{C}_{i}\}_{i=1}^{n}$, $K$,
$\{\alpha_{k}\}_{k=1}^{K}$, $\eta\geq1$, $\Delta$, and $\epsilon$.

\textbf{Output}: $\{\hat{r}_{i}^{(K+1)}\}_{i=1}^{n}$.

Initialize $\{r_{i}^{(1)}\}_{i=1}^{n}$ arbitrarily. Let $\hat{r}_{i}^{(1)}=r_{i}^{(1)}$
for all $i\in[n]$ and $\theta_{k}=(\eta+1)/(\eta+k)$ for $k\in[K]$.

For $k=1,2,\dots,K$, repeat:
\begin{enumerate}
\item If $k=1$, then set $w_{k}=0$; else draw a random vector $w_{k}\in\mathbb{R}^{T}$
from the distribution (proportional to) $\exp\left(-\frac{2\epsilon\norm{w_{k}}}{K(K-1)L\Delta}\right)$
\item Compute~$\hat{p}^{(k)}:=\nabla U\left({\textstyle \sum_{i=1}^{n}r_{i}^{(k)}}\right)+w_{k}$.
\item For $i\in[n]$, compute: 
\begin{align*}
r_{i}^{(k+1)} & :=\Pi_{\mathcal{C}_{i}}(r_{i}^{(k)}-\alpha_{k}\hat{p}^{(k)}),\\
\hat{r}_{i}^{(k+1)} & :=(1-\theta_{k})\hat{r}_{i}^{(k)}+\theta_{k}r_{i}^{(k+1)}.
\end{align*}

\end{enumerate}
\caption{Differentially private distributed projected gradient descent.\label{alg:DP_ev_charging}}
\end{algorithm}

In the first half of the paper, we present the main algorithmic result
of this paper, a differentially private distributed algorithm (Algorithm~\ref{alg:DP_ev_charging})
for solving the constrained optimization problem~(\ref{eq:dco_prob}).
The constant $\Delta$ that appears in the input of Algorithm~\ref{alg:DP_ev_charging}
is defined as
\begin{multline}
\Delta:=\max_{i\in[n]}\max\Bigl\{\norm{\Pi_{\mathcal{C}_{i}}(r)-\Pi_{\mathcal{C}_{i}'}(r)}\colon\\
r\in\mathbb{R}^{T},\ \mathcal{C}_{i}\text{ and }\mathcal{C}_{i}'\text{ satisfy }\rho_{\mathcal{C}}(\mathcal{C}_{i},\mathcal{C}_{i}')\leq\delta\mathcal{C}\Bigr\}.\label{eq:def_delta}
\end{multline}
In other words, $\Delta$ can be viewed as a bound on the global $\ell_{2}$-sensitivity
of the projection operator $\Pi_{\mathcal{C}_{i}}$ to changes in
$\mathcal{C}_{i}$ for all $i\in[n]$. Later, we will illustrate how
to compute $\Delta$ using the case of EV charging. 

Compared to the (non-private) distributed algorithm (Algorithm~\ref{alg:Distributed-EV-charging}),
the key difference in Algorithm~\ref{alg:DP_ev_charging} is the
introduction of random perturbations in the gradients (step 2) that
convert $p^{(k)}$ into a noisy gradient $\hat{p}^{(k)}$. The noisy
gradients $(\hat{p}^{(1)},\hat{p}^{(2)},\dots,\hat{p}^{(K)})$ can
be viewed as a randomized mechanism~$M_{p}$ that approximates the
original gradients $p=(p^{(1)},p^{(2)},\dots,p^{(K)})$. In Section~\ref{sec:A-differentially-private},
we will prove that the noisy gradients (as a mechanism) $M_{p}:=(\hat{p}^{(1)},\hat{p}^{(2)},\dots,\hat{p}^{(K)})$
\emph{preserve $\epsilon$-differential privacy} and hence solve Problem~\ref{prob:Find-a-randomized}. 
\begin{thm}
\label{thm:Algorithm-ensures-that}Algorithm~\ref{alg:DP_ev_charging}
ensures that $M_{p}:=(\hat{p}^{(1)},\hat{p}^{(2)},\dots,\hat{p}^{(K)})$
preserves $\epsilon$-differential privacy under the adjacency relation
given by Definition~\ref{def:adj_dco}.
\end{thm}
Algorithm~\ref{alg:DP_ev_charging} can be viewed as an instance
of stochastic gradient descent that terminates after $K$ iterations.
We will henceforth refer to Algorithm~\ref{alg:DP_ev_charging} as
\emph{differentially private distributed projected gradient descent}.
The step size $\alpha_{k}$ is chosen as $\alpha_{k}=c/\sqrt{k}$
for some $c>0$. The purpose of the additional variables $\{\hat{r}_{i}^{(k)}\}_{k=1}^{K}$
is to implement the polynomial-decay averaging method in order to
improve the convergence rate, which is a common practice in stochastic
gradient descent~\cite{shamir2013stochastic}; introducing $\{\hat{r}_{i}^{(k)}\}_{k=1}^{K}$
does not affect privacy. The parameter $\eta\geq1$ is used for controlling
the averaging weight $\theta_{k}$. Details on choosing $\eta$ can
be found in Shamir and Zhang~\cite{shamir2013stochastic}.

Like most iterative optimization algorithms, stochastic gradient descent
only converges in a probabilistic sense as the number of iterations
$K\to\infty$. In practice, the number of iterations is always finite,
so that it is desirable to analyze the suboptimality for a finite~$K$.
In Section~\ref{sec:A-differentially-private}, we provide an analysis
on the expected suboptimality of Algorithm~\ref{alg:DP_ev_charging}.
\begin{thm}
\label{thm:subopt_bound}The expected suboptimality of Algorithm~\ref{alg:DP_ev_charging}
after $K$ iterations is bounded as follows:
\begin{multline}
\mathbb{E}\left[U\left({\textstyle \sum_{i=1}^{n}\hat{r}_{i}^{(K+1)}}\right)-U^{*}\right]\\
\qquad\leq\mathcal{O}\left(\eta\sqrt{n}\rho\left(\frac{G}{\sqrt{K}}+\frac{\sqrt{2}TK^{3/2}L\Delta}{2\epsilon}\right)\right),\label{eq:ev_subopt}
\end{multline}
where~$U^{*}$ is the optimal value of problem~(\ref{eq:dco_prob}),
and 
\begin{align*}
\rho & =\max\left\{ \sqrt{{\textstyle \sum_{i=1}^{n}}\norm{r_{i}}^{2}}\colon r_{i}\in\mathcal{C}_{i},\ i\in[n]\right\} ,\\
G & =\max\{\Vert\nabla U\left({\textstyle \sum_{i=1}^{n}r_{i}}\right)\Vert\colon r_{i}\in\mathcal{C}_{i},\ i\in[n]\}.
\end{align*}

\end{thm}

\subsection{Results for the case of EV charging}

Having presented and analyzed the algorithm for a general distributed
constrained optimization problem, in the second half of the paper,
we illustrate how Algorithm~\ref{alg:DP_ev_charging} can be applied
to the case of EV charging. In Section~\ref{sec:Sensitivity-computation},
we demonstrate how to compute $\Delta$ using the case of EV charging
as an example. %
We show that $\Delta$ can be bounded by $\delta r$ and $\delta E$
that appear in~(\ref{eq:adj_ev}) as described below.
\begin{thm}
\label{thm:The-global--sensitivity}Consider the example of EV charging
(as described in Section~\ref{sub:Example-application:-EV}). For
any $i\in[n]$, the global $\ell_{2}$-sensitivity of the projection
operator $\Pi_{\mathcal{C}_{i}(\bar{r}_{i},E_{i})}$ with respect
to changes in $(\bar{r}_{i},E_{i})$ is bounded by
\[
\Delta\leq2\delta r+\delta E,
\]
where $\delta r$ and $\delta E$ are specified in the adjacency relation
given by~(\ref{eq:adj_ev}). 
\end{thm}
The suboptimality analysis given in Theorem~\ref{thm:subopt_bound}
can be further refined in the case of EV charging. The special form
of $U$ given by~(\ref{eq:cost_quad}) allows obtaining an upper
bound on suboptimality as given in Corollary~\ref{cor:quad_cost}
below.
\begin{cor}
\label{cor:quad_cost}For the cost function $U$ given by~(\ref{eq:cost_quad}),
the expected suboptimality of Algorithm~\ref{alg:DP_ev_charging}
is bounded as follows:
\begin{equation}
\mathbb{E}\left[U\left({\textstyle \sum_{i=1}^{n}\hat{r}_{i}^{(K+1)}}\right)-U^{*}\right]\leq\mathcal{O}\left(\eta T{}^{1/4}(\Delta/n\epsilon)^{1/4}\right).\label{eq:ev_subopt_quad}
\end{equation}

\end{cor}
This upper bound shows the trade-off between privacy and performance.
As $\epsilon$ decreases, more privacy is preserved but at the expense
of increased suboptimality. On the other hand, this increase in suboptimality
can be mitigated by introducing more participating users (i.e., by
increasing $n$), which coincides with the common intuition that it
is easier to achieve privacy as the number of users $n$ increases.

\section{Differentially private distributed projected gradient descent\label{sec:A-differentially-private}}

In this section, we give the proof that the modified distributed projected
gradient descent algorithm (Algorithm~\ref{alg:DP_ev_charging})
preserves $\epsilon$-differential privacy. %
In the proof, we will extensively use results from differential privacy
such as the Laplace mechanism and the adaptive sequential composition
theorem.

\subsection{Review: Results from differential privacy}

The introduction of additive noise in step 2 of Algorithm~\ref{alg:DP_ev_charging}
is based on a variant of the widely used Laplace mechanism in differential
privacy. The Laplace mechanism operates by introducing additive noise
according to the $\ell_{p}$-sensitivity ($p\geq1$) of a numerical
query~$q\colon\mathcal{D}\to\mathbb{R}^{m}$ (for some dimension
$m$), which is defined as follows.
\begin{defn}
[$\ell_p$-sensitivity]For any query $q\colon\mathcal{D}\to\mathbb{R}^{m}$,
the $\ell_{p}$\emph{-sensitivty} of $q$ under the adjacency relation
$\mathrm{Adj}$ is defined as
\begin{align*}
 & \Delta_{q}:=\max\{\norm{q(D)-q(D')}_{p}\colon\\
 & \qquad D,D'\in\mathcal{D}\text{ s.t. }\mathrm{Adj}(D,D')\}.
\end{align*}

\end{defn}
\noindent Note that the $\ell_{p}$-sensitivity of $q$ does not depend
on a specific database~$D$. In this paper, we will use the Laplace
mechanism for bounded $\ell_{2}$-sensitivity. 
\begin{prop}
[Laplace mechanism~\cite{dwork2006calibrating}]\label{prop:laplace}Consider
a query $q\colon\mathcal{D}\to\mathbb{R}^{m}$ whose $\ell_{2}$-sensitivity
is $\Delta_{q}$. Define the mechanism~$M_{q}$ as $M_{q}(D):=q(D)+w$,
where $w$ is an $m$-dimensional random vector whose probability
density function is given by $p_{w}(w)\propto\exp(-\epsilon\norm{w}/\Delta_{q})$.
Then, the mechanism $M_{q}$ preserves $\epsilon$-differential privacy.
\end{prop}
As a basic building block in differential privacy, the Laplace mechanism
allows construction of the differentially private distributed projected
gradient descent algorithm described in Algorithm~\ref{alg:DP_ev_charging}
through \emph{adaptive sequential composition}.
\begin{prop}
[Adaptive seqential composition~\cite{dwork2013algorithmic}]\label{prop:seq_composition_adp}Consider
a sequence of mechanisms $\{M_{k}\}_{k=1}^{K}$, in which the output
of $M_{k}$ may depend on $M_{1},M_{2},\dots,M_{k-1}$ as described
below:
\[
M_{k}(D)=M_{k}(D,M_{1}(D),M_{2}(D),\dots,M_{k-1}(D)).
\]
Suppose~$M_{k}(\cdot,a_{1},a_{2},\dots,a_{k-1})$ preserves $\epsilon_{k}$-differential
privacy for any $a_{1}\in\mathrm{range}(M_{1}),\dots,a_{k-1}\in\mathrm{range}(M_{k-1}).$
Then, the $K$-tuple mechanism $M:=(M_{1},M_{2},\dots,M_{K})$ preserves
$\epsilon$-differential privacy for $\epsilon=\sum_{k=1}^{K}\epsilon_{k}$.
\end{prop}

\subsection{Proof on that Algorithm~\ref{alg:DP_ev_charging} preserves $\epsilon$-differential
privacy }

Using the adaptive sequential composition theorem, we can show that
Algorithm~\ref{alg:DP_ev_charging} preserves $\epsilon$-differential
privacy. We can view the $K$-tuple mechanism $M_{p}:=(\hat{p}^{(1)},\hat{p}^{(2)},\dots,\hat{p}^{(K)})$
as a sequence of mechanisms $\{\hat{p}^{(k)}\}_{k=1}^{K}$. The key
is to compute the $\ell_{2}$-sensitivity of $p^{(k)}:=\nabla U\left({\textstyle \sum_{i=1}^{n}r_{i}^{(k)}}\right)$,
denoted by $\Delta^{(k)}$, when the outputs of $\hat{p}^{(1)},\hat{p}^{(2)},\dots,\hat{p}^{(k-1)}$
are given, so that we can obtain a differentially private mechanism
$\hat{p}^{(k)}$ by applying the Laplace mechanism on $p^{(k)}$ according
to $\Delta^{(k)}$.
\begin{lem}
\label{lem:When-the-outputs}In Algorithm~\ref{alg:DP_ev_charging},
when the outputs of $\hat{p}^{(1)},\hat{p}^{(2)},\dots,\hat{p}^{(k-1)}$
are given, the $\ell_{2}$-sensitivity of $p^{(k)}:=\nabla U\left({\textstyle \sum_{i=1}^{n}r_{i}^{(k)}}\right)$
satisfies $\Delta^{(k)}=(k-1)L\Delta$.\end{lem}
\begin{IEEEproof}
See Appendix~\ref{sub:Proof-of-Lemma-Sensivity}.
\end{IEEEproof}
With Lemma~\ref{lem:When-the-outputs} at hand, we now show that
Algorithm~\ref{alg:DP_ev_charging} preserves $\epsilon$-differential
privacy (Theorem~\ref{thm:Algorithm-ensures-that}, Section~\ref{sec:overview_of_results}). 
\begin{IEEEproof}
(of Theorem~\ref{thm:Algorithm-ensures-that}) For any $k\in[K]$,
when the outputs of $\hat{p}^{(1)},\hat{p}^{(2)},\dots,\hat{p}^{(k-1)}$
are given, we know from Proposition~\ref{prop:laplace} that $\hat{p}^{(k)}$
preserves $\epsilon_{k}$-differential privacy, where $\epsilon_{k}$
satisfies $\epsilon_{1}=0$ and for $k>1$, 
\[
\epsilon_{k}/\Delta^{(k)}=\frac{2\epsilon}{K(K-1)L\Delta}.
\]
Use the expression of $\Delta^{(k)}$ from Lemma~\ref{lem:When-the-outputs}
to obtain 
\[
\epsilon_{k}=\frac{2(k-1)\epsilon}{K(K-1)}.
\]
Using the adaptive sequential composition theorem, we know that the
privacy of $M_{p}:=(\hat{p}^{(1)},\hat{p}^{(2)},\dots,\hat{p}^{(K)})$
is given by $\sum_{k=1}^{K}\epsilon_{k}=\frac{2\epsilon}{K(K-1)}\sum_{k=1}^{K}(k-1)=\epsilon$,
which completes the proof.
\end{IEEEproof}

\subsection{Suboptimality analysis: Privacy-performance trade-off}

As a consequence of preserving privacy, we only have access to noisy
gradients $\{\hat{p}^{(k)}\}_{k=1}^{K}$ rather than the exact gradients
$\{p^{(k)}\}_{k=1}^{K}$. Recall that the additive noise $w_{k}$
in step 2 of Algorithm~\ref{alg:DP_ev_charging} has zero mean. In
other words, the noisy gradient $\hat{p}^{(k)}$ is an unbiased estimate
of $p^{(k)}$, which allows us to view Algorithm~\ref{alg:DP_ev_charging}
as an instantiation of stochastic gradient descent. As we mentioned
in Section~\ref{sec:overview_of_results}, it is in fact a variant
of stochastic gradient descent that uses polynomial-decay averaging
for better convergence. The stochastic gradient descent method (with
polynomial-decay averaging), which is described in Algorithm~\ref{alg:SGD_polyave},
can be used for solving the following optimization problem:
\[
\underset{x}{\optmin}\quad f(x)\qquad\optst\quad x\in\mathcal{X},
\]
where $x\in\mathbb{R}^{m}$ and $\mathcal{X}\subset\mathbb{R}^{m}$
for certain dimensions~$m$. Proposition~\ref{prop:SGD_subopt}
(due to Shamir and Zhang~\cite{shamir2013stochastic}) gives an upper
bound of the expected suboptimality after finitely many steps for
the stochastic gradient descent algorithm (Algorithm~\ref{alg:SGD_polyave}).

\begin{algorithm}
\textbf{Input}: $f$, $\mathcal{X}$, $K$, $\{\alpha_{k}\}_{k=1}^{K}$,
and $\eta\ge1$.

\textbf{Output}: $\hat{x}^{(K+1)}$.

Initialize $x^{(1)}$ and $k=1$. Let $\hat{x}^{(1)}=x^{(1)}$ and
$\theta_{k}=(\eta+1)/(\eta+k)$ for $k\in[K]$.

For $k=1,2,\dots,K$, repeat:
\begin{enumerate}
\item Compute an unbiased subgradient~$\hat{g}_{k}$ of~$f$ at $x^{(k)}$,
i.e., $\mathbb{E}[\hat{g}_{k}]\in\partial f(x^{(k)})$.
\item Update $x^{(k+1)}:=\Pi_{\mathcal{X}}(x^{(k)}-\alpha_{k}\hat{g}_{k})$
and $\hat{x}^{(k+1)}:=(1-\theta_{k})\hat{x}^{(k)}+\theta_{k}x^{(k+1)}$.
\end{enumerate}
\caption{Stochastic gradient descent with polynomial-decay averaging.\label{alg:SGD_polyave}}
\end{algorithm}

\begin{prop}
[Shamir and Zhang~\cite{shamir2013stochastic}]\label{prop:SGD_subopt}Suppose
$\mathcal{X}\subset\mathbb{R}^{m}$ is a convex set and $f\colon\mathbb{R}^{m}\to\mathbb{R}$
is a convex function. Assume that there exist~$\rho$ and $\widehat{G}$
such that $\sup_{x,x'\in\mathcal{X}}\norm{x-x'}\leq\rho$ and $\max_{1\leq k\leq K}\mathbb{E}\norm{\hat{g}_{k}}^{2}\leq\widehat{G}^{2}$
for $\{\hat{g}_{k}\}_{k=1}^{K}$ given by step 1 of Algorithm~\ref{alg:SGD_polyave}.
If the step sizes are chosen as $\alpha_{k}=c/\sqrt{k}$ for some~$c>0$
, then for any $K>1$, it holds that
\begin{equation}
\mathbb{E}(f(\hat{x}^{(K+1)})-f^{*})\leq\mathcal{O}\left(\frac{\eta(\rho^{2}/c+c\widehat{G}^{2})}{\sqrt{K}}\right),\label{eq:generic_sgd_subopt}
\end{equation}
where~$f^{*}=\inf_{x\in\mathcal{X}}f(x)$. 
\end{prop}
A tighter upper bound can be obtained from~(\ref{eq:generic_sgd_subopt})
by optimizing the right-hand side of~(\ref{eq:generic_sgd_subopt})
over the constant~$c$.
\begin{cor}
\label{cor:SGD_subopt} Under the same setting as Proposition~\ref{prop:SGD_subopt},
the suboptimality bound for Algorithm~\ref{alg:SGD_polyave} is given
by
\begin{equation}
\mathbb{E}(f(\hat{x}^{(K+1)})-f^{*})\leq\mathcal{O}\left(\eta\frac{\rho\widehat{G}}{\sqrt{K}}\right),\label{eq:generic_sgd_subopt-1}
\end{equation}
which is achieved by choosing $c=\rho/\widehat{G}$. 
\end{cor}
By applying Corollary~\ref{cor:SGD_subopt}, we are able obtain the
bound of suboptimality for Algorithm~\ref{alg:DP_ev_charging} as
given by Theorem~\ref{thm:subopt_bound}.
\begin{IEEEproof}
(of Theorem~\ref{thm:subopt_bound}) In order to apply Corollary~\ref{cor:SGD_subopt},
we need to compute $\rho$ and $\widehat{G}$ for Algorithm~\ref{alg:DP_ev_charging}.
The constant $\rho$ can be obtained as 
\[
\rho=\max\left\{ \sqrt{{\textstyle \sum_{i=1}^{n}}\norm{r_{i}}^{2}}\colon r_{i}\in\mathcal{C}_{i},\ i\in[n]\right\} .
\]
Recall that the definition of $\widehat{G}$ is given by $\widehat{G}^{2}:=\max_{k}\mathbb{E}\norm{\hat{g}_{k}}^{2}$.
It can be verified that the gradient $\hat{g}_{k}$ of the objective
function~$U\left({\textstyle \sum_{i=1}^{n}r_{i}^{(k)}}\right)$
with respective to $(r_{1}^{(k)},r_{2}^{(k)},\dots,r_{n}^{(k)})$
is given by $\hat{g}_{k}=[\hat{p}^{(k)},\hat{p}^{(k)},\dots,\hat{p}^{(k)}]$,
which is formed by repeating $\hat{p}^{(k)}$ for $n$ times, so that
we have $\widehat{G}^{2}=n\cdot\max_{k}\mathbb{E}\norm{\hat{p}^{(k)}}^{2}$.
Using the expression of $\hat{p}^{(k)}$, we have 
\begin{align*}
\widehat{G} & =\sqrt{n}\cdot\max_{k\in[K]}\sqrt{\norm{p^{(k)}}^{2}+\mathbb{E}\norm{w_{k}}^{2}}\\
 & \leq\sqrt{n}\cdot\max_{k\in[K]}\left\{ \Vert p^{(k)}\Vert+\sqrt{\mathbb{E}\norm{w_{k}}^{2}}\right\} \\
 & \leq\sqrt{n}(G+\sqrt{2}TK^{2}L\Delta/2\epsilon),
\end{align*}
where in the last step we have used the fact that
\begin{align*}
\mathbb{E}\norm{w_{k}}^{2} & =\mathrm{var}\norm{w_{k}}^{2}+\left(\mathbb{E}\norm{w_{k}}\right)^{2}\\
 & =T(\Delta^{(k)}/\epsilon_{k})^{2}+T^{2}(\Delta^{(k)}/\epsilon_{k})^{2}\\
 & \leq2T^{2}(\Delta^{(k)}/\epsilon_{k})^{2}
\end{align*}
and
\[
\Delta^{(k)}/\epsilon_{k}=\frac{K(K-1)L\Delta}{2\epsilon}\leq\frac{K^{2}L\Delta}{2\epsilon}.
\]
Substitute the expression of $\widehat{G}$ into~(\ref{eq:generic_sgd_subopt-1})
to obtain the result. 
\end{IEEEproof}
As~$K$ increases, the first term in~(\ref{eq:ev_subopt}) decreases,
whereas the second term in~(\ref{eq:ev_subopt}) increases. It is
then foreseeable that there exists an optimal choice of~$K$ that
minimizes the expected suboptimality. 
\begin{cor}
\label{cor:The-minimum-expected}The expected suboptimality of Algorithm~\ref{alg:DP_ev_charging}
is bounded as follows:
\begin{multline}
\mathbb{E}\left[U\left({\textstyle \sum_{i=1}^{n}\hat{r}_{i}^{(K+1)}}\right)-U^{*}\right]\\
\leq\mathcal{O}\left(\eta T{}^{1/4}n^{1/2}\rho(G^{3}L\Delta/\epsilon)^{1/4}\right),\label{eq:ev_subopt_final}
\end{multline}
where $U^{*}$, $\rho$, and $G$ are given by Theorem~\ref{thm:subopt_bound}.
The bound~(\ref{eq:ev_subopt}) is achieved by choosing $K=(\sqrt{2}G\epsilon/3TL\Delta)^{1/2}$.\end{cor}
\begin{IEEEproof}
The result can be obtained by optimizing the right-hand side of~(\ref{eq:ev_subopt})
over~$K$.
\end{IEEEproof}
However, since it is generally impossible to obtain a tight bound
for~$\rho$ and~$\widehat{G}$, optimizing~$K$ according to Corollary~\ref{cor:The-minimum-expected}
usually does not give the best~$K$ in practice; numerical simulation
is often needed in order to find the best~$K$ for a given problem.
We will demonstrate how to choose $K$ optimally later using numerical
simulations in Section~\ref{sec:Numerical-simulations}.

\section{Sensitivity computation: The case of EV charging\label{sec:Sensitivity-computation}}

So far, we have shown that Algorithm~\ref{alg:DP_ev_charging} (specifically,
the mechanism~$M_{p}$ consisting of the gradients $(\hat{p}^{(1)},\hat{p}^{(2)},\dots,\hat{p}^{(K)})$
that are broadcast to every participating user) preserves $\epsilon$-differential
privacy. The magnitude of the noise~$w_{k}$ introduced to the gradients
depends on $\Delta$, which is the sensitivity of the projection operator
$\Pi_{\mathcal{C}_{i}}$ as defined in~(\ref{eq:def_delta}). In
order to implement Algorithm~\ref{alg:DP_ev_charging}, we need to
compute~$\Delta$ explicitly. In the next, we will illustrate how
to compute $\Delta$ using the case of EV charging as an example.
We will give an expression for $\Delta$ that depends on the constants
$\delta r$ and $\delta E$ appearing in the adjacency relation~(\ref{eq:adj_ev}).
Since $\delta r$ and $\delta E$ are part of the privacy requirement,
one can choose $\Delta$ accordingly once the privacy requirement
has been determined.

\subsection{Overview}

The input of Algorithm~\ref{alg:DP_ev_charging} includes the constant~$\Delta$
as described by~(\ref{eq:def_delta}), which bounds the global $\ell_{2}$-sensitivity
of the projection operator $\Pi_{\mathcal{C}_{i}(\bar{r}_{i},E_{i})}$
with respect to changes in $(\bar{r}_{i},E_{i})$. In this section,
we will derive an explicit expression of $\Delta$ for the case of
EV charging. Using tools in sensitivity analysis of optimization problems,
we are able to establish the relationship between~$\Delta$ and the
constants $\delta r$ and $\delta E$ that appear in the adjacency
relation~(\ref{eq:adj_ev}) used in EV charging.

Recall that for any $r\in\mathbb{R}^{T}$, the output of the projection
operation $\Pi_{\mathcal{C}_{i}(\bar{r}_{i},E_{i})}(r)$ is the optimal
solution to the constrained least-squares problem
\begin{alignat}{2}
 & \underset{r_{i}}{\optmin}\quad &  & \frac{1}{2}\Vert r_{i}-r\Vert^{2}\label{eq:ev_charging_prob_alt-1}\\
 & \optst\quad &  & 0\preceq r_{i}\preceq\bar{r}_{i},\qquad\mathbf{1}^{T}r_{i}=E_{i}.\nonumber 
\end{alignat}
Define the $\ell_{2}$-sensitivity for a fixed $r$ as 
\begin{multline*}
\Delta_{r}:=\max_{i\in[n]}\max\Bigl\{\norm{\Pi_{\mathcal{C}_{i}(\bar{r}_{i}.E_{i})}(r)-\Pi_{\mathcal{C}_{i}(\bar{r}_{i}',E_{i}')}(r)}\colon\\
(\bar{r}_{i},E_{i})\text{ and }(\bar{r}_{i}',E_{i}')\text{ satisfy }\eqref{eq:adj_ev}\Bigr\};
\end{multline*}
it can be verified that $\Delta=\max_{r\in\mathbb{R}^{T}}\Delta_{r}$.
In the following, we will establish the relationship between $\Delta_{r}$
and $(\delta r,\delta E)$; we will also show that $\Delta_{r}$ does
not depend on the choice of $r$, so that $\Delta=\Delta_{r}$ for
any $r\in\mathbb{R}^{T}$. For notational convenience, we consider
the following least-squares problem:
\begin{alignat}{2}
 & \underset{x}{\optmin}\quad &  & \frac{1}{2}\Vert x-x_{0}\Vert^{2}\label{eq:ev_charging_prob_alt}\\
 & \optst\quad &  & 0\preceq x\preceq a,\qquad\mathbf{1}^{T}x=b,\nonumber 
\end{alignat}
where $x_{0}$, $a$, and $b$ are given constants. If we let 
\[
x_{0}=r_{0},\qquad a=\bar{r}_{i},\qquad b=E_{i},
\]
then problem~(\ref{eq:ev_charging_prob_alt}) is mapped back to problem~(\ref{eq:ev_charging_prob_alt-1}).
We also have $\mathbf{1}^{T}a\geq b$ based on the assumption as described
in~(\ref{eq:feasiblity}). Denote the optimal solution of problem~(\ref{eq:ev_charging_prob_alt})
by $x^{*}(a,b)$. Since our purpose is to derive an expression for
$\Delta_{r}$ when $r$ is fixed, we also treat $x_{0}$ as fixed
and has dropped the dependence of $x^{*}$ on $x_{0}$. Our goal is
to bound the global solution sensitivity with respect to changes in
$a$ and $b$, i.e., 
\begin{equation}
\norm{x^{*}(a',b')-x^{*}(a,b)},\label{eq:global_sensitivity_ab}
\end{equation}
for any $(a,b)$ and $(a',b')$. We will proceed by first bounding
the \emph{local solution sensitivities} $\partial_{a}x^{*}$ and $\partial_{b}x^{*}$
with respect to $a$ and $b$. Then, we will obtain a bound on the
global solution, sensitivity~(\ref{eq:global_sensitivity_ab}) through
integration of $\partial_{a}x^{*}$ and $\partial_{b}x^{*}$.

\subsection{Local solution sensitivity of nonlinear optimization problems}

We begin by reviewing existing results on computing local solution
sensitivity of nonlinear optimization problems. Consider a generic
nonlinear optimization problem parametrized by~$\theta\in\mathbb{R}$
described as follows:
\begin{alignat}{2}
 & \underset{x\in\mathbb{R}^{n}}{\optmin}\quad &  & f(x;\theta)\label{eq:generic_nonlinear}\\
 & \optst\quad &  & g_{i}(x;\theta)\leq0,\quad i\in[p]\nonumber \\
 &  &  & h_{j}(x;\theta)=0,\quad j\in[q],\nonumber 
\end{alignat}
whose Lagrangian can be expressed as 
\[
L(x,\lambda,\nu;\theta)=f(x;\theta)+\sum_{i=1}^{p}\lambda_{i}g_{i}(x;\theta)+\sum_{j=1}^{q}\nu_{j}h_{j}(x;\theta),
\]
where $\lambda$ and $\nu$ are the Lagrange multipliers associated
with constraints $\{g_{i}\}_{i=1}^{p}$ and $\{h_{j}\}_{j=1}^{q}$,
respectively. If there exists a set~$\Theta\subset\mathbb{R}$ such
that the optimal solution is unique for all $\theta\in\Theta$, then
the optimal solution of problem~(\ref{eq:generic_nonlinear}) can
be defined as a function~$x^{*}\colon\Theta\to\mathbb{R}^{n}$. This
condition on the uniqueness of optimal solution holds for problem~(\ref{eq:ev_charging_prob_alt}),
since the objective function therein is strictly convex. 

Denote by $\lambda^{*}$ and $\nu^{*}$ the optimal Lagrange multipliers.
Under certain conditions described in Proposition~\ref{thm:opt_sensitivity}
below, the partial derivatives $\partial_{\theta}x^{*}$, $\partial_{\theta}\lambda^{*}$,
and $\partial_{\theta}\nu^{*}$ exist; these partial derivatives will
also be referred to as \emph{local solution sensitivity} of problem~(\ref{eq:generic_nonlinear}).
\begin{prop}
[Fiacco~\cite{fiacco1976sensitivity}]\label{thm:opt_sensitivity}Let
$(x^{*},\lambda^{*},\nu^{*})$ be the primal-dual optimal solution
of problem~(\ref{eq:generic_nonlinear}). Suppose the following conditions
hold.
\begin{enumerate}
\item $x^{*}$ is a locally unique optimal primal solution.
\item The functions~$ $$f$, $\{g_{i}\}_{i=1}^{p},$ and $\{h_{j}\}_{j=1}^{q}$
are twice continuously differentiable in~$x$ and differentiable
in~$\theta$. 
\item The gradients $\{\nabla g_{i}(x^{*})\colon g_{i}(x^{*})=0,\: i\in[p]\}$
of the active constraints and the gradients $\{\nabla h_{j}(x^{*})\colon j\in[q]\}$
are linearly independent.
\item Strict complementary slackness holds: $\lambda_{i}^{*}>0$ when $g_{i}(x^{*},\theta)=0$
for all $i\in[p]$.
\end{enumerate}

Then the local sensitivity $(\partial_{\theta}x^{*},\partial_{\theta}\lambda^{*},\partial_{\theta}\nu^{*})$
of problem~(\ref{eq:generic_nonlinear}) exists and is continuous
in a neighborhood of $\theta$. Moreover, $(\partial_{\theta}x^{*},\partial_{\theta}\lambda^{*},\partial_{\theta}\nu^{*})$
is uniquely determined by the following: 
\[
\nabla^{2}L\cdot\partial_{\theta}x^{*}+\sum_{i=1}^{p}\nabla g_{i}\cdot\partial_{\theta}\lambda_{i}^{*}+\sum_{j=1}^{q}\nabla h_{j}\cdot\partial_{\theta}\nu_{j}^{*}+\partial_{\theta}(\nabla L)=0
\]
and
\begin{align*}
\lambda_{i}\nabla g_{i}\cdot\partial_{\theta}x^{*}+g_{i}\partial_{\theta}\lambda_{i}^{*}+\lambda_{i}^{*}\partial_{\theta}g_{i} & =0,\quad i\in[p]\\
\nabla h_{j}\cdot\partial_{\theta}x^{*}+\partial_{\theta}h_{j} & =0,\quad j\in[q].
\end{align*}

\end{prop}

\subsection{Solution sensitivity of the distributed EV charging problem}

We begin by computing the local solution sensitivities $\partial_{a}x^{*}$
and $\partial_{b}x^{*}$ for problem~(\ref{eq:ev_charging_prob_alt})
using Proposition~\ref{thm:opt_sensitivity}. After $\partial_{a}x^{*}$
and $\partial_{b}x^{*}$ are obtained, the global solution sensitivity~(\ref{eq:global_sensitivity_ab})
can be obtained through integration of the local sensitivity. For
convenience, we compute the global solution sensitivity in $a$ and
$b$ separately and combine the results in the end in the proof of
Theorem~\ref{thm:The-global--sensitivity}. 

One major difficulty in applying Proposition~\ref{thm:opt_sensitivity}
is that it requires strict complementary slackness, which unfortunately
does not hold for all values of $a$ and $b$. We will proceed by
deriving the local solution sensitivities assuming that strict complementary
slackness holds. Later, we will show that strict complementary slackness
only fails at finitely many locations on the integration path, so
that the integral remains unaffected. 

We shall proceed by computing the solution sensitivity for $a$ and
$b$ separately. First of all, we assume that $a$ is fixed and solve
for the global solution sensitivity of $x^{*}$ in $b$, defined as
\begin{equation}
\norm{x^{*}(a,b')-x^{*}(a,b)}.\label{eq:global_sensitivity_b}
\end{equation}
When strict complementary slackness holds, the following lemma gives
properties of the local solution sensitivity of problem~(\ref{eq:ev_charging_prob_alt})
with respect to $b$. 
\begin{lem}
[Local solution sensitivity in $b$]\label{lem:When-strict-complementary}When
 strict complementary slackness holds, the local solution sensitivity
$\partial_{b}x^{*}$ of problem~(\ref{eq:ev_charging_prob_alt})
satisfies 
\[
\partial_{b}x^{*}\succeq0\qquad\text{and}\qquad\mathbf{1}^{T}\partial_{b}x^{*}=1.
\]
\end{lem}
\begin{IEEEproof}
See Appendix~\ref{sub:appendix-b1}.
\end{IEEEproof}
The following lemma shows that the condition is only violated for
a finite number of values of~$b$, so that it will still be possible
to obtain the global sensitivity~(\ref{eq:global_sensitivity_b})
through integration.
\begin{lem}
\label{lem:riemann_integrable}The set of possible values of $b$
in problem~(\ref{eq:ev_charging_prob_alt}) for which strict complementary
slackness fails to hold is finite.\end{lem}
\begin{IEEEproof}
See Appendix~\ref{sub:appendix-b2}.
\end{IEEEproof}
The implication of Lemma~\ref{lem:riemann_integrable} is that the
local solution sensitivity~$\partial_{b}x^{*}$ exists everywhere
except at finitely many location. It is also possible to show that
the optimal solution $x^{*}(a,b)$ is continuous in $b$ (cf.~Berge~\cite[page 116]{berge1963topological},
Dantzig et al.~\cite{dantzig1967continuity}, or Tropp et al.~\cite[Appendix I]{tropp2005designing});
the continuity of $x^{*}(a,b)$ in $b$ and together with Lemma~\ref{lem:riemann_integrable}
imply that $\partial_{b}x^{*}$ is Riemann integrable so that we can
obtain the global solution sensitivity through integration.
\begin{prop}
[Global solution sensivitity in $b$]\label{thm:global_solution_sensitivity}For
any $a$, $b$ and $b'$ that satisfy $\mathbf{1}^{T}a\geq b$ and
$\mathbf{1}^{T}a\geq b'$, we have
\[
\norm{x^{*}(a,b')-x^{*}(a,b)}_{1}=|b'-b|.
\]
\end{prop}
\begin{IEEEproof}
Without loss of generality, assume~$b'>b$. Since $x^{*}(a,b)$ is
continuous in~$b$, and the partial derivative $\partial_{b}x^{*}$
exists except at finitely many points according to Lemma~\ref{lem:riemann_integrable},
we know that $\partial_{b}x_{i}^{*}$ is Riemann integrable for all
$i\in[T]$, so that 
\[
x_{i}^{*}(a,b')-x_{i}^{*}(a,b)=\int_{b}^{b'}\partial_{b}x_{i}^{*}(a,b)\, db
\]
according to the fundamental theorem of calculus. Using the fact $\partial_{b}x_{i}^{*}\geq0$
as given by Lemma~\ref{lem:When-strict-complementary}, we have
\[
x_{i}^{*}(a,b')-x_{i}^{*}(a,b)\geq0.
\]
Then, we have
\begin{multline*}
\norm{x^{*}(a,b')-x^{*}(a,b)}_{1}=\sum_{i=1}^{T}|x_{i}^{*}(a,b')-x_{i}^{*}(a,b)|\\
=\sum_{i=1}^{T}(x_{i}^{*}(a,b')-x_{i}^{*}(a,b))=\int_{b}^{b'}\mathbf{1}^{T}\partial_{b}x^{*}(a,b)\, db.
\end{multline*}
Using the fact $\mathbf{1}^{T}\partial_{b}x^{*}=1$ given by Lemma~\ref{lem:When-strict-complementary},
we obtain
\[
\norm{x^{*}(a,b')-x^{*}(a,b)}_{1}=b'-b=|b'-b|.
\]

\end{IEEEproof}
Having computed the solution sensitivity in $b$, we now assume that
$b$ is fixed and solve for the global solution sensitivity of $x^{*}$
in $a$, defined as
\begin{equation}
\norm{x^{*}(a',b)-x^{*}(a,b)}.\label{eq:global_sensitivity_a}
\end{equation}
When strict complementary slackness holds, the following lemma gives
properties of the local solution sensitivity of problem~(\ref{eq:ev_charging_prob_alt})
with respect to $a$. 
\begin{lem}
[Local solution sensivitity in $a$]\label{lem:When-strict-complementary-1}When
strict complementary slackness holds, the local solution sensitivity
$\partial_{a}x^{*}$ of problem~(\ref{eq:ev_charging_prob_alt})
satisfies
\[
\sum_{i=1}^{T}\norm{\partial_{a_{i}}x^{*}}_{1}\leq2.
\]
 \end{lem}
\begin{IEEEproof}
See Appendix~\ref{sub:appendix-b3}.
\end{IEEEproof}
Similar to computing the sensitivity in $b$, we can obtain the global
solution sensitivity in~$a$ by integration of the local sensitivity.
For convenience, we choose the integration path~$L$ from any $a$
to $a'$ such that only one component of $a$ is varied at a time.
Namely, the path~$L$ is given by 
\begin{align}
 & L\colon(a_{1},a_{2},\dots,a_{T})\to(a_{1}',a_{2},\dots,a_{T})\to\cdots\nonumber \\
 & \qquad\to(a_{1}',a_{2}',\dots,a_{T}').\label{eq:int_path}
\end{align}
For convenience, we also define the subpaths $L_{1},L_{2},\dots,L_{T}$
such that 
\begin{equation}
L_{i}\colon(a_{1}',\dots,a_{i-1}',a_{i},\dots,a_{T})\to(a_{1}',\dots,a_{i-1}',a_{i}',\dots,a_{T}).\label{eq:int_subpath}
\end{equation}
It is also possible to establish the fact that $\partial_{a}x^{*}$
exists excepts for a finite number of locations along the integration
path~$L$. Note that we do not need to check whether strict complementary
slackness holds for the constraint $x_{i}\geq0$, since in this case
$\partial_{a_{i}}x^{*}$ always exists (in fact, $\partial_{a_{i}}x^{*}=0$);
instead, we only need to check strict complementary slackness 
\begin{equation}
\mu_{i}^{*}>0\quad\text{when}\quad x_{i}^{*}=a_{i}\quad\text{for all }i\label{eq:strict_cs_a}
\end{equation}
associated with the constraint $x\preceq a$.
\begin{lem}
\label{lem:riemann_integrable-1}When constrained on the integration
path~$L$ given by~(\ref{eq:int_path}), the set of possible values
of $a$ in problem~(\ref{eq:ev_charging_prob_alt}) for which the
strict complementary slackness condition~(\ref{eq:strict_cs_a})
fails to hold is finite.\end{lem}
\begin{IEEEproof}
See Appendix~\ref{sub:appendix-b4}.
\end{IEEEproof}
Lemma~\ref{lem:riemann_integrable-1} guarantees that $\partial_{a}x^{*}$
is Riemann integrable along~$L$, so that we can obtain the global
solution sensitivity~(\ref{eq:global_sensitivity_a}) through integration.
\begin{prop}
[Global solution sensivitity in $a$]\label{prop:global_solution_sensitivity-1}For
any given $a$, $a'$, and $b$ that satisfy $\mathbf{1}^{T}a\geq b$
and $\mathbf{1}^{T}a'\geq b$, we have\textup{
\[
\norm{x^{*}(a',b)-x^{*}(a,b)}_{1}\leq2\norm{a'-a}_{1}.
\]
}\end{prop}
\begin{IEEEproof}
Similar to the proof of Proposition~\ref{thm:global_solution_sensitivity},
we can show that $\partial_{a}x_{i}^{*}$ is Riemann integrable using
both Lemma~\ref{lem:riemann_integrable-1} and the fact that $x^{*}$
is continuous in~$a$. Then, we can define
\[
I_{ij}:=\int_{L_{j}}\partial_{a}x_{i}^{*}(a,b)\cdot d\ell,
\]
which is the line integral of the vector field $\partial_{a_{j}}x_{i}^{*}(\cdot,b)$
along the path $L_{j}$. Define $x_{ij}^{a}:=\partial_{a_{j}}x_{i}^{*}(\cdot,b)$.
Using the definition of $L_{j}$ as given by~(\ref{eq:int_subpath}),
we can write $I_{ij}$ as 
\[
I_{ij}=\int_{a_{j}}^{a_{j}'}x_{ij}^{a}(a_{1}',a_{2}',\dots,a_{j-1}',a_{j},\dots,a_{T})\, da_{j}.
\]
Then, we have $ $
\[
x_{i}^{*}(a',b)-x_{i}^{*}(a,b)=\int_{L}\partial_{a}x_{i}^{*}(a,b)\cdot d\ell=\sum_{j=1}^{T}I_{ij}
\]
and consequently
\[
|x_{i}^{*}(a',b)-x_{i}^{*}(a,b)|\leq\sum_{j=1}^{T}|I_{ij}|.
\]
Substituting the expression of $|x_{i}^{*}(a',b)-x_{i}^{*}(a,b)|$
into the the global sensitivity expression~(\ref{eq:global_sensitivity_a}),
we obtain
\begin{align}
 & \norm{x^{*}(a',b)-x^{*}(a,b)}_{1}=\sum_{i=1}^{T}|x_{i}^{*}(a',b)-x_{i}^{*}(a,b)|\nonumber \\
 & \qquad\leq\sum_{i=1}^{T}\sum_{j=1}^{T}|I_{ij}|=\sum_{j=1}^{T}\sum_{i=1}^{T}|I_{ij}|.\label{eq:sum_Iij}
\end{align}
Note that we have 
\[
|I_{ij}|\leq\int_{\underline{a}_{j}}^{\bar{a}_{j}}\left|x_{ij}^{a}(a_{1}',\dots,a_{j-1}',a_{j},\dots,a_{T})\right|\, da_{j},
\]
where $\underline{a}_{j}:=\min(a_{j},a_{j}')$ and $\bar{a}_{j}:=\max(a_{j},a_{j}')$,
so that 
\begin{align*}
\sum_{i=1}^{T}|I_{ij}| & \leq\sum_{i=1}^{T}\int_{\underline{a}_{j}}^{\bar{a}_{j}}\left|x_{ij}^{a}(a_{1}',\dots,a_{j-1}',a_{j},\dots,a_{T})\right|\, da_{j}\\
 & =\int_{\underline{a}_{j}}^{\bar{a}_{j}}\sum_{i=1}^{T}\left|x_{ij}^{a}(a_{1}',\dots,a_{j-1}',a_{j},\dots,a_{T})\right|\, da_{j}\\
 & =\int_{\underline{a}_{j}}^{\bar{a}_{j}}\norm{\partial_{a_{j}}x^{*}((a_{1}',\dots,a_{j-1}',a_{j},\dots,a_{T}),b)}_{1}\, da_{j}\\
 & :=\bar{I}_{j}.
\end{align*}
Using Lemma~\ref{lem:When-strict-complementary-1}, we can show that
$\bar{I}_{j}$ satisfies
\[
\bar{I}_{j}\le\int_{\underline{a}_{j}}^{\bar{a}_{j}}2\, da_{j}=2|a_{j}'-a_{j}|.
\]
Substitute the above into~(\ref{eq:sum_Iij}) to obtain 
\[
\norm{x^{*}(a',b)-x^{*}(a,b)}_{1}\leq2\sum_{j=1}^{T}|a_{j}'-a_{j}|=2\norm{a'-a}_{1},
\]
which completes the proof.
\end{IEEEproof}
We are now ready to prove Theorem~\ref{thm:The-global--sensitivity}
using results from Propositions~\ref{thm:global_solution_sensitivity}
and~\ref{prop:global_solution_sensitivity-1}.
\begin{IEEEproof}
(of Theorem~\ref{thm:The-global--sensitivity}) Consider problem~(\ref{eq:ev_charging_prob_alt}).
By combining Propositions~\ref{thm:global_solution_sensitivity}
and~\ref{prop:global_solution_sensitivity-1}, we can obtain the
global solution sensitivity with respect to both $a$ and $b$ as
defined by~(\ref{eq:global_sensitivity_ab}). Consider any given
$(a,b)$ and $(a',b')$ that satisfy $\mathbf{1}^{T}a\geq b$ and
$\mathbf{1}^{T}a'\geq b'$. Without loss of generality, we assume
that $\mathbf{1}^{T}a'\geq\mathbf{1}^{T}a$, so that $\mathbf{1}^{T}a'\geq b$;
this implies that the corresponding optimization problem(\ref{eq:ev_charging_prob_alt})
is feasible, and the optimal solution $x^{*}(a',b)$ is well-defined.
Then, we have 
\begin{align}
 & \norm{x^{*}(a',b')-x^{*}(a,b)}_{1}\nonumber \\
 & \qquad=\norm{x^{*}(a',b')-x^{*}(a',b)+x^{*}(a',b)-x^{*}(a,b)}_{1}\nonumber \\
 & \qquad\leq\norm{x^{*}(a',b')-x^{*}(a',b)}_{1}+\norm{x^{*}(a',b)-x^{*}(a,b)}_{1}\nonumber \\
 & \qquad\leq\norm{b'-b}_{1}+2\norm{a'-a}_{1}.\label{eq:joint_sensitivity_ab}
\end{align}

By letting $x_{0}=r_{0}$, $a=\bar{r}_{i}$, and $b=E_{i}$, we can
map problem~(\ref{eq:ev_charging_prob_alt}) back to problem~(\ref{eq:ev_charging_prob_alt-1}).
Recall that $x^{*}$ is defined as the optimal solution of problem~(\ref{eq:ev_charging_prob_alt})
for a given $x_{0}$. Then, the inequality~(\ref{eq:joint_sensitivity_ab})
implies that for any given $r$, the $\ell_{2}$-sensitivity
\[
\Delta_{r}\leq2\norm{\bar{r}_{i}'-\bar{r}_{i}}_{1}+\norm{E_{i}'-E_{i}}_{1}=2\delta r+\delta E.
\]
However, since the right-hand side of the above inequality does not
depend on $r$, we have 
\[
\Delta=\max_{r\in\mathbb{R}^{T}}\Delta_{r}\leq2\delta r+\delta E,
\]
which completes the proof.\end{IEEEproof}
\begin{rem}
Alternatively, one may use the fact $\Pi_{\mathcal{C}_{i}}(\cdot)\in\mathcal{C}_{i}$
to obtain a bound on $\Delta$. Recall that for any $r_{i}\in\mathcal{C}_{i}$
we have 
\[
\norm{r_{i}}\leq\norm{\bar{r}_{i}}\quad\text{and}\quad\norm{r_{i}}\leq\norm{\hat{r}}_{1}=E_{i}.
\]
Then, we have 
\[
\Delta\leq\max_{i\in[n]}\left(\norm{\bar{r}_{i}}+\norm{\bar{r}_{i}'}\right)\leq2\max_{i\in[n]}\norm{\bar{r}_{i}}+\delta r
\]
and 
\[
\Delta\leq\max_{i\in[n]}(E_{i}+E_{i}')\leq2\max_{i\in[n]}E_{i}+\delta E.
\]
However, this bound can be quite loose in practice. Since the magnitude
of~$w_{k}$ in Algorithm~\ref{alg:DP_ev_charging} is proportional
to $\Delta$, a loose bound on $\Delta$ implies introducing more
noise to the gradient $p^{(k)}$ than what is necessary for preserving
$\epsilon$-differential privacy. As we have already seen in Section~(\ref{sub:Suboptimality-analysis-EV}),
the noise magnitude is closely related to the performance loss of
Algorithm~\ref{alg:DP_ev_charging} caused by preserving privacy,
and less noise is always desired for minimizing such a loss.
\end{rem}

\subsection{Revisited: Suboptimality analysis\label{sub:Suboptimality-analysis-EV}}

The bound~(\ref{eq:ev_subopt_final}) given by Corollary~\ref{cor:The-minimum-expected}
does not clearly indicate the dependence of suboptimality on the number
of participating users $n$, because $\rho$, $G$, and $L$ also
depend on $n$. In order to reveal the dependence on $n$, we will
further refine the suboptimality bound~(\ref{eq:ev_subopt_final})
for the specific objective function~$U$ given in~(\ref{eq:cost_quad}).
The resulting suboptimality bound is shown in Corollary~\ref{cor:quad_cost}.
\begin{IEEEproof}
(of Corollary~\ref{cor:quad_cost}) Define $r_{\max}=\max_{i\in[n]}\norm{r_{i}}.$
Then, we have $\rho\leq\sqrt{n}r_{\max}$. For $U$ given by~(\ref{eq:cost_quad}),
its gradient can be computed as
\[
\nabla U\left({\textstyle \sum_{i=1}^{n}r_{i}}\right)=\frac{1}{m}\left(d+{\textstyle \sum_{i=1}^{n}r_{i}}/m\right),
\]
so that
\[
\norm{\nabla U\left({\textstyle \sum_{i=1}^{n}r_{i}}\right)}\leq\frac{1}{m}\left(\norm{d}+{\textstyle \sum_{i=1}^{n}\norm{r_{i}}}/m\right).
\]
Then, we obtain
\begin{align*}
G & :=\max\{\Vert\nabla U\left({\textstyle \sum_{i=1}^{n}r_{i}}\right)\Vert\colon r_{i}\in\mathcal{C}_{i},\: i\in[n]\}.\\
 & \leq\frac{1}{m}\left(\norm{d}+n{\textstyle r_{\max}}/m\right)\\
 & =\frac{\gamma}{n}\left(\norm{d}+\gamma{\textstyle r_{\max}}\right).
\end{align*}
In order to compute the Lipschitz constant $L$ for $\nabla U$, note
that for any $x,y\in\mathbb{R}^{T}$, we have
\[
\norm{\nabla U(x)-\nabla U(y)}=\frac{1}{m}\norm{x/m-y/m}=\frac{1}{m^{2}}\norm{x-y},
\]
so that we obtain $L=1/m^{2}=\gamma^{2}/n^{2}$. Substitute $\rho$,
$G$, and $L$ into~(\ref{eq:ev_subopt_final}) to obtain~(\ref{eq:ev_subopt_quad}).
Note that we have dropped the dependence on $d$, $r_{\max}$, and
$\gamma$ for brevity, since we are most interested in the relationship
between suboptimality and $(n,\epsilon)$.
\end{IEEEproof}
The suboptimality bound~(\ref{eq:ev_subopt_quad}) indicates how
performance (cost) is affected by incorporating privacy. As $\epsilon$
decreases, the level of privacy is elevated but at the expense of
sacrificing performance as a result of increased suboptimality. This
increase in suboptimality can be mitigated by introducing more participating
users (i.e., by increasing $n$) as predicted by the bound~(\ref{eq:ev_subopt_quad});
this coincides with the common intuition that it is easier to achieve
privacy as the number of users $n$ increases when only aggregate
user information is available. Indeed, in the distributed EV charging
algorithm, the gradients $p^{(k)}$ is a function of the aggregate
load profile $\sum_{i=1}^{n}r_{i}^{(k)}$.

\section{Numerical simulations\label{sec:Numerical-simulations}}

We consider the cost function as given by~(\ref{eq:cost_quad}).
The base load $d$ is chosen according to the data provided in Gan
et al.~\cite{gan2013optimal}. The scheduling horizon is divided
into $52$ time slots of $15$ minutes. We consider a large pool of
EVs ($n=100,000$) in a large residential area ($m=500,000$). For
computational efficiency, instead of assigning a different charging
specification $(\bar{r}_{i},E_{i})$ to every user~$i\in[n]$, we
divide the users into $N$ ($N\ll n$) groups and assign the same
charging specification for every user in the same group. If we choose
the same initial conditions $r_{i}^{(1)}$ for all users in the same
group, the projected gradient descent update (step 3, Algorithm~\ref{alg:DP_ev_charging})
becomes identical for all users in the group, so that the projection~$\Pi_{\mathcal{C}_{i}}$
only needs to be computed once for a given group. We choose $N=100$
and draw $(\bar{r}_{j},E_{j})$ for all $j\in[N]$ as follows. The
entries of $\bar{r}_{j}$ are drawn independently from a Bernoulli
distribution, where $\bar{r}_{j}(t)=3.3$~kW with probability $0.5$
and $\bar{r}_{j}(t)=0$~kW with probability $0.5$. The amount of
energy $E_{j}$ is drawn from the uniform distribution on the interval
$[28,40]$~(kW). Note that $E_{j}$ has been normalized against~$\Delta T=0.25$~h
to match the unit of $\bar{r}_{j}$; in terms of energy required,
this implies that each vehicle needs an amount between $[28,40]\:\text{kW}\times0.25\:\text{h}=[7,10]\:\text{kWh}$
by the end of the scheduling horizon. 

The constants $\delta r$ and $\delta E$ in the adjacency relation~(\ref{eq:adj_ev})
are determined as follows. We choose $\delta r=3.3\times4=13.2$~kW,
so that the privacy of any events spanning less than $4$ time slots
(i.e., 1 hour) can be preserved; we choose $\delta E=40-28=12$~kW,
which corresponds to the maximum difference in $E_{j}$ ($j\in[N]$)
as $E_{j}$ varies between the interval~$[28,40]$ (kW). The other
parameters in Algorithm~\ref{alg:DP_ev_charging} are chosen as follows:
$L=1/m^{2}$ (as computed in Section~\ref{sub:Suboptimality-analysis-EV}),
$\eta=1$, whereas~$\epsilon$, $K$, and $c$ will vary among different
numerical experiments.

\begin{figure}
\begin{centering}
\includegraphics[width=2in]{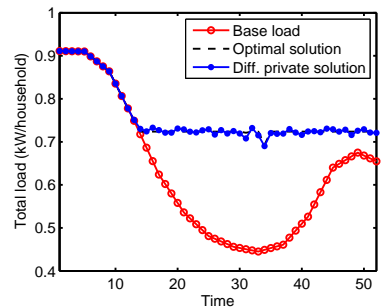}
\par\end{centering}

\caption{A typical output of the differentially private distributed EV charging
algorithm (Algorithm~\ref{alg:DP_ev_charging}) with $\epsilon=0.1$
compared to the optimal solution of problem~(\ref{eq:dco_prob}).
The other parameters used in the simulation are: $K=6$ and $c=10$.}

\label{fig:A-typical-output}
\end{figure}

\begin{figure}
\begin{centering}
\includegraphics[width=2in]{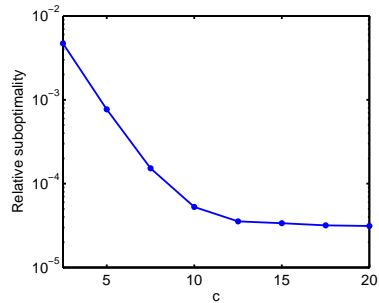}
\par\end{centering}

\caption{Relative suboptimality of the differentially private distributed EV
charging algorithm (Algorithm~\ref{alg:DP_ev_charging}) as a function
of the step size constant~$c$. The other parameters used in the
simulations are: $\epsilon=0.1$ and~$K=6$.\label{fig:Suboptimality-vs-c}}
\end{figure}

Fig.~\ref{fig:A-typical-output} plots a typical output from Algorithm~\ref{alg:DP_ev_charging}
with $\epsilon=0.1$, alongside the optimal solution of problem~(\ref{eq:dco_prob}).
The dip at $t=34$ that appears both in the differentially private
solution and the optimal solution is due to the constraint imposed
by $\bar{r}_{i}(t)$. Because of the noise introduced in the gradients,
the differentially private solution given by Algorithm~\ref{alg:DP_ev_charging}
exhibits some additional fluctuations compared to the optimal solution.
The constant $c$ that determines the step sizes is found to be insensitive
for $c\in[10,20]$ as shown in Fig.~\ref{fig:Suboptimality-vs-c},
so that we choose $c=10$ in all subsequent simulations. In Fig.~\ref{fig:Suboptimality-vs-c},
the relative suboptimality is defined as $\left[U\left({\textstyle \sum_{i=1}^{n}\hat{r}_{i}^{(K+1)}}\right)-U^{*}\right]/U^{*}$,
which is obtained by normalizing the suboptimality against $U^{*}$.

\begin{figure}
\begin{centering}
\includegraphics[width=2in]{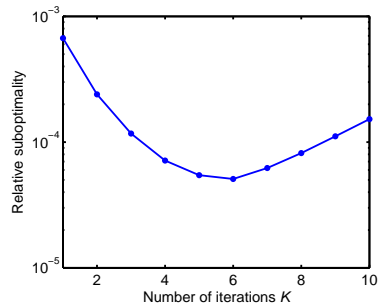}
\par\end{centering}

\caption{Relative suboptimality of the differentially private distributed EV
charging algorithm (Algorithm~\ref{alg:DP_ev_charging}) as a function
of the number of iterations~$K$. The other parameters used in the
simulations are: $\epsilon=0.1$ and~$c=10$.\label{fig:Suboptimality-vs-K}}
\end{figure}

Fig.~\ref{fig:Suboptimality-vs-K} shows the relative suboptimality
as a function of~$K$. It can be seen from Fig.~\ref{fig:Suboptimality-vs-K}
that an optimal choice of~$K$ exists, which coincides with the result
of Theorem~\ref{thm:subopt_bound}.

\begin{figure}
\begin{centering}
\includegraphics[width=2in]{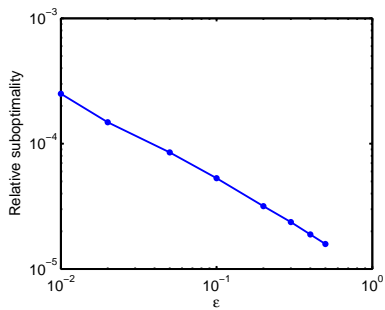}
\par\end{centering}

\caption{Relative suboptimality of Algorithm~\ref{alg:DP_ev_charging} as
a function of~$\epsilon$. Larger~$\epsilon$ implies that less
privacy is preserved. The slope is approximately $-0.698$ (compared
to the theoretical bound $-0.25$ as given by Corollary~\ref{cor:quad_cost}).
All simulations use $c=10$. The number of iterations~$K$ is optimized
for every choice of~$\epsilon$. \label{fig:Suboptimality-vs-eps}}
\end{figure}

Fig.~\ref{fig:Suboptimality-vs-eps} shows the dependence of the
relative suboptimality on~$\epsilon$. A separate experiment for
investigating the dependence on~$n$ is not performed, since changing~$n$
is expected to have a similar effect as changing~$\epsilon$ according
to Corollary~\ref{cor:quad_cost}. As the privacy requirement becomes
less stringent (i.e., as~$\epsilon$ grows), the suboptimality of
Algorithm~\ref{alg:DP_ev_charging} improves, which coincides qualitatively
with the bound given in Corollary~\ref{cor:quad_cost}. One can quantify
the relationship between the suboptimality and $\epsilon$ from the
slope of the curve in Fig.~\ref{fig:Suboptimality-vs-eps}; if the
slope is $s$, then the relationship between the suboptimality and
$\epsilon$ is given by~$\mathcal{O}(\epsilon^{s})$. By performing
linear regression on the curve in Fig.~\ref{fig:Suboptimality-vs-eps},
we can obtain the slope as $s\approx-0.698$, which is in contrast
to $-0.25$ given by Corollary~\ref{cor:quad_cost}. This implies
that the suboptimality of Algorithm~\ref{alg:DP_ev_charging} decreases
faster than the rate given by Corollary~\ref{cor:quad_cost} as $\epsilon$
increases. In other words, the bound given by Corollary~\ref{cor:quad_cost}
is likely to be loose; this is possibly due to the fact that the result
on the suboptimality of stochastic gradient descent (Proposition~\ref{prop:SGD_subopt})
does not consider additional properties of the objective function
such as strong convexity. 

\section{Conclusions}

This paper develops an $\epsilon$-differentially private algorithm
for distributed constrained optimization. The algorithm preserves
privacy in the specifications of user constraints by adding noise
to the public coordination signal (i.e., gradients). By using the
sequential adaptive composition theorem, we show that the noise magnitude
is determined by the sensitivity of the projection operation~$\Pi_{\mathcal{C}}$,
where $\mathcal{C}$ is the parameterized set describing the user
constraints. By viewing the projection operation as a least-squares
problem, we are able to compute the sensitivity of $\Pi_{\mathcal{C}}$
through a solution sensitivity analysis of optimization problems.
We demonstrate how this sensitivity can be computed in the case of
EV charging. 

We also analyze the trade-off between privacy and performance of the
algorithm through results on suboptimality analysis of the stochastic
gradient descent method. Specifically, in the case of EV charging,
the expected suboptimality of the $\epsilon$-differentially private
optimization algorithm with~$n$ participating users is upper bounded
by~$\mathcal{O}((n\epsilon)^{-1/4})$. For achieving the best suboptimality,
both the suboptimality analysis and numerical simulations show that
there exists an optimal choice for the number iterations: too few
iterations affects the convergence behavior, whereas too many iterations
leads to too much noise in the gradients. Simulations have indicated
that the bound $\mathcal{O}((n\epsilon)^{-1/4})$ is likely not tight.
One future direction is to derive a tighter bound for similar distributed
optimization problems using information-theoretic approaches (e.g.,~\cite{duchi2013local,raginsky2011information}).

\bibliographystyle{abbrv}
\bibliography{ref,/Users/hanshuo/git/privacy/ref,ref2,ref_games}

\appendices

\section{Proof of Lemma~\ref{lem:When-the-outputs}\label{sub:Proof-of-Lemma-Sensivity}}

Consider any adjacent $D$ and $D'$ such that $\mathcal{C}_{j}=\mathcal{C}_{j}'$
for all $j\neq i$. We will first show that when the outputs of $\hat{p}^{(1)},\hat{p}^{(2)},\dots,\hat{p}^{(k-1)}$
are given, we have
\begin{align*}
\norm{r_{i}^{(k)}(D')-r_{i}^{(k)}(D)} & \leq(k-1)\Delta,\\
\norm{r_{j}^{(k)}(D')-r_{j}^{(k)}(D)} & =0,\quad\forall j\neq i.
\end{align*}
We will prove the above result by induction. For $k=1$, we have $\norm{r_{i}^{(1)}(D')-r_{i}^{(1)}(D)}=0$
for all $i\in[n]$. 

Consider the case when $k>1$. For notational convenience, we define
for $i\in[n]$, 
\begin{equation}
v_{i}^{(k-1)}(D):=r_{i}^{(k-1)}(D)-\alpha_{k-1}\hat{p}^{(k-1)}.\label{eq:v_k}
\end{equation}
In~(\ref{eq:v_k}), we have used the fact that the output of $\hat{p}^{(k-1)}$
is given so that the dependence of $\hat{p}^{(k-1)}$ on $D$ is dropped
according to the adaptive sequential composition theorem (Proposition~\ref{prop:seq_composition_adp}).
Then, for all $j\neq i$, we have
\begin{align*}
 & \norm{r_{j}^{(k)}(D')-r_{j}^{(k)}(D)}\\
 & \qquad=\norm{\Pi_{\mathcal{C}_{j}'}(v_{j}^{(k-1)}(D'))-\Pi_{\mathcal{C}_{j}}(v_{j}^{(k-1)}(D))}\\
 & \qquad=\norm{\Pi_{\mathcal{C}_{j}}(v_{j}^{(k-1)}(D'))-\Pi_{\mathcal{C}_{j}}(v_{j}^{(k-1)}(D))}\\
 & \qquad\leq\norm{v_{j}^{(k-1)}(D')-v_{j}^{(k-1)}(D)}\\
 & \qquad=\norm{r_{j}^{(k-1)}(D')-r_{j}^{(k-1)}(D)}=0
\end{align*}
and
\begin{align*}
 & \norm{r_{i}^{(k)}(D')-r_{i}^{(k)}(D)}\\
 & \qquad=\norm{\Pi_{\mathcal{C}_{i}'}(v_{i}^{(k-1)}(D'))-\Pi_{\mathcal{C}_{i}}(v_{i}^{k-1)}(D))}\\
 & \qquad\leq\norm{\Pi_{\mathcal{C}_{i}'}(v_{i}^{(k-1)}(D'))-\Pi_{\mathcal{C}_{i}}(v_{i}^{k-1)}(D'))}\\
 & \qquad\qquad+\norm{\Pi_{\mathcal{C}_{i}}(v_{i}^{(k-1)}(D'))-\Pi_{\mathcal{C}_{i}}(v_{i}^{k-1)}(D))}\\
 & \qquad\leq\Delta+\norm{v_{i}^{(k-1)}(D')-v_{i}^{(k-1)}(D)}\\
 & \qquad=\Delta+\norm{r_{i}^{(k-1)}(D')-r_{i}^{(k-1)}(D)}\leq(k-1)\Delta,
\end{align*}
where we have used the induction hypothesis
\begin{align*}
\norm{r_{i}^{(k-1)}(D')-r_{i}^{(k-1)}(D)} & \leq(k-2)\Delta,\\
\norm{r_{j}^{(k-1)}(D')-r_{j}^{(k-1)}(D)} & =0,\quad\forall j\neq i.
\end{align*}
Then, the $\ell_{2}$-sensitivity of $p^{(k)}$ can be computed as
follows:
\begin{multline*}
\norm{p^{(k)}(D')-p^{(k)}(D)}\\
\leq L\norm{\sum_{i=1}^{n}\left[r_{i}^{(k)}(D')-r_{i}^{(k)}(D)\right]}\leq L(k-1)\Delta.
\end{multline*}
Since the above results hold for all $i$ such that $D$ and $D'$
satisfy $\mathcal{C}_{j}=\mathcal{C}_{j}'$ for all $j\neq i$, we
have 
\begin{align*}
\Delta^{(k)} & :=\max_{D,D'\colon\mathrm{Adj}(D,D')}\norm{p^{(k)}(D')-p^{(k)}(D)}\\
 & =\max_{i\in[n]}\max\Bigl\{\norm{p^{(k)}(D')-p^{(k)}(D)}\colon D,D'\text{ satisfy}\\
 & \qquad\mathcal{C}_{j}=\mathcal{C}_{j}'\text{ for all }j\neq i\Bigr\}\\
 & =L(k-1)\Delta.
\end{align*}

\section{Proofs on the local solution sensitivities\label{sec:Proofs-on-the-local}}

\subsection{Proof of Lemma~\ref{lem:When-strict-complementary}\label{sub:appendix-b1}}

The Lagrangian of problem~(\ref{eq:ev_charging_prob_alt}) can be
written as
\begin{equation}
L(x,\lambda,\mu,\nu)=\frac{1}{2}\norm{x-x_{0}}^{2}-\lambda^{T}x+\mu^{T}(x-a)+\nu(b-\mathbf{1}^{T}x).\label{eq:proj_lagrangian}
\end{equation}
Denote by $\lambda^{*}$, $\mu^{*}$, and $\nu^{*}$ the corresponding
optimal Lagrange multipliers. It can be verified that all conditions
in Proposition~\ref{thm:opt_sensitivity} hold. Apply Proposition~\ref{thm:opt_sensitivity}
to obtain
\begin{align}
\partial_{b}x^{*}-\partial_{b}\lambda^{*}+\partial_{b}\mu^{*}-\partial_{b}\nu^{*}\cdot\mathbf{1} & =0\label{eq:kkt_sensitivity_ev1}\\
\mathbf{1}^{T}\partial_{b}x^{*} & =1\label{eq:kkt_sensitivity_ev2}\\
\lambda_{i}^{*}\cdot\partial_{b}x_{i}^{*}+x_{i}^{*}\cdot\partial_{b}\lambda_{i}^{*} & =0,\quad i\in[T]\label{eq:kkt_sensitivity_ev3}\\
\mu_{i}^{*}\cdot\partial_{b}x_{i}^{*}+(x_{i}^{*}-a_{i})\cdot\partial_{b}\mu_{i}^{*} & =0,\quad i\in[T].\label{eq:kkt_sensitivity_ev4}
\end{align}
Strict complementary slackness implies that either (1) $x_{i}^{*}=0$
and $\lambda_{i}^{*}>0$, so that $\partial_{b}x_{i}^{*}=0$ according
to~(\ref{eq:kkt_sensitivity_ev3}); or (2) $x_{i}^{*}\neq0$ and
$\lambda_{i}^{*}=0$, so that $\partial_{b}\lambda_{i}^{*}=0$ also
according to~(\ref{eq:kkt_sensitivity_ev3}). In other words, under
strict complementary slackness, condition~(\ref{eq:kkt_sensitivity_ev3})
is equivalent to 
\begin{equation}
\partial_{b}\lambda_{i}^{*}\cdot\partial_{b}x_{i}^{*}=0.\label{eq:kkt_sensitivity_ev3-2}
\end{equation}
Similarly, we can rewrite condition~(\ref{eq:kkt_sensitivity_ev4})
as
\begin{equation}
\partial_{b}\mu_{i}^{*}\cdot\partial_{b}x_{i}^{*}=0.\label{eq:kkt_sensitivity_ev4-2}
\end{equation}
Conditions~(\ref{eq:kkt_sensitivity_ev3-2}) and~(\ref{eq:kkt_sensitivity_ev4-2})
imply that one and only one of the following is true for any $i\in[T]$:
(1) $\partial_{b}x_{i}^{*}=0$; (2) $\partial_{b}\lambda_{i}^{*}=0$
and $\partial_{b}\mu_{i}^{*}=0$. Define~$\mathcal{I}:=\{i\colon\partial_{b}x_{i}^{*}\neq0\}$,
and we have 
\begin{equation}
\sum_{i\in\mathcal{I}}\partial_{b}x_{i}^{*}=1\label{eq:dxdb_sum}
\end{equation}
from~(\ref{eq:kkt_sensitivity_ev2}). Note that~(\ref{eq:kkt_sensitivity_ev1})
implies that for all $i,j\in[T]$, 
\[
\partial_{b}x_{i}^{*}-\partial_{b}\lambda_{i}^{*}+\partial_{b}\mu_{i}^{*}=\partial_{b}x_{j}^{*}-\partial_{b}\lambda_{j}^{*}+\partial_{b}\mu_{j}^{*}.
\]
Since $\partial_{b}\lambda_{i}^{*}=0$ and $\partial_{b}\mu_{i}^{*}=0$
for all $i\in\mathcal{I}$, we have $\partial_{b}x_{i}^{*}=\partial_{b}x_{j}^{*}$
for all $i,j\in\mathcal{I}$ and hence $\partial_{b}x_{i}^{*}=1/|\mathcal{I}|$
for all $i\in\mathcal{I}$ according to~(\ref{eq:dxdb_sum}). On
the other hand, from the definition of $\mathcal{I}$, we have $\partial_{b}x_{i}^{*}=0$
for all $i\notin\mathcal{I}$. In summary, we have $\partial_{b}x^{*}\succeq0$,
which complete the proof.

\subsection{Proof of Lemma~\ref{lem:riemann_integrable}\label{sub:appendix-b2}}

The optimality conditions for problem~(\ref{eq:ev_charging_prob_alt})
imply that
\begin{align}
x^{*}-\lambda^{*}+\mu^{*}-\nu^{*}\mathbf{1} & =x_{0}\label{eq:kkt_ev1}\\
\mathbf{1}^{T}x^{*} & =b\label{eq:kkt_ev2}\\
\lambda_{i}^{*}x_{i}^{*} & =0,\quad i\in[T]\label{eq:kkt_ev3}\\
\mu_{i}^{*}(x_{i}^{*}-a_{i}) & =0,\quad i\in[T].\label{eq:kkt_ev4}
\end{align}
Suppose strict complementary slackness fails for a certain value of~$b$.
Denote the set of indices of the constraints that violate strict complementary
slackness by $\mathcal{I}_{\lambda}:=\{i\colon\lambda_{i}^{*}=0,\: x_{i}^{*}=0\}$
and $\mathcal{I}_{\mu}:=\{i\colon\mu_{i}^{*}=0,\: x_{i}^{*}=a_{i}\}$.
If both $\mathcal{I}_{\lambda}$ and $\mathcal{I}_{\mu}$ are empty,
then strict complementary slackness holds for all constraints. 
\begin{IEEEproof}
When $\mathcal{I}_{\lambda}$ is non-empty, we know from~(\ref{eq:kkt_ev4})
that $\mu_{i}^{*}=0$ for all $i\in\mathcal{I}_{\lambda}$. For any~$i\in\mathcal{I}_{\lambda}$,
substitute $x_{i}^{*}=0$, $\lambda_{i}^{*}=0$, and $\mu_{i}^{*}=0$
into~(\ref{eq:kkt_ev1}) to obtain $\nu^{*}=x_{0,i}$%
. For any other $j\notin\mathcal{I}_{\lambda}$, one of the following
three cases must hold: (1) $x_{j}^{*}=0$; (2) $x_{j}^{*}=a_{j}$;
(3) $0<x_{j}^{*}<a_{j}$. Consider a partition $(\mathcal{I}_{1},\mathcal{I}_{2},\mathcal{I}_{3})$
of the set $[n]\backslash\mathcal{I}_{\lambda}$ as follows: 
\begin{gather*}
\mathcal{I}_{1}:=\{j\colon x_{j}^{*}=0\},\qquad\mathcal{I}_{2}:=\{j\colon x_{j}^{*}=a_{j}\},\\
\mathcal{I}_{3}:=\{j\colon0<x_{j}^{*}<a_{j}\}.
\end{gather*}
For any $j\in\mathcal{I}_{3}$, we have $\lambda_{j}^{*}=\mu_{j}^{*}=0$
from~(\ref{eq:kkt_sensitivity_ev3}) and~(\ref{eq:kkt_sensitivity_ev4}),
so that we have $x_{j}^{*}=\nu^{*}+x_{0,j}$ according to~(\ref{eq:kkt_sensitivity_ev1}).
Then, we can write using~(\ref{eq:kkt_ev2}) 
\begin{align}
b & =\mathbf{1}^{T}x^{*}\nonumber \\
 & =\sum_{i\in\mathcal{I_{\lambda}}}x_{i}^{*}+\sum_{j\in\mathcal{I}_{1}}x{}_{j}^{*}+\sum_{j\in\mathcal{I}_{2}}x_{j}^{*}+\sum_{j\in\mathcal{I}_{3}}x_{j}^{*}\nonumber \\
 & =|\mathcal{I}_{\lambda}|\nu^{*}+0+\sum_{j\in\mathcal{I}_{2}}a_{j}+\sum_{j\in\mathcal{I}_{3}}(\nu^{*}+x_{0,j}).\label{eq:b_partition}
\end{align}
Since both $a$ and $x_{0}$ are fixed, we know that the choice of
$\nu^{*}=x_{0,i}$ (for any $i\in\mathcal{I}_{\lambda}$) is finite.
By enumerating all finitely many partitions $(\mathcal{I}_{\lambda},\mathcal{I}_{1},\mathcal{I}_{2},\mathcal{I}_{3})$
of $[T]$, we know that $b$ can only take finitely many values according
to~(\ref{eq:b_partition}). The proof is similar for the case when
$\mathcal{I}_{\mu}$ is nonempty by making use of~(\ref{eq:kkt_ev3}).
When both $\mathcal{I}_{\lambda}$ and $\mathcal{I}_{\mu}$ are nonempty,
the possible values of $b$ are given by the intersection of those
when only one of $\mathcal{I}_{\lambda}$ and $\mathcal{I}_{\mu}$
is empty; hence the number of possible values is also finite.
\end{IEEEproof}

\subsection{Proof of Lemma~\ref{lem:When-strict-complementary-1}\label{sub:appendix-b3}}

Similar to the proof of Lemma~\ref{lem:When-strict-complementary},
we apply Proposition~\ref{thm:opt_sensitivity} using the Lagrangian
as given by~(\ref{eq:proj_lagrangian}). We can show that the following
holds for all $i\in[T]$:
\begin{align}
\partial_{a_{i}}x_{i}^{*}-\partial_{a_{i}}\lambda_{i}^{*}+\partial_{a_{i}}\mu_{i}^{*}-\partial_{a_{i}}\nu^{*} & =0\label{eq:kkt_ev_a1}\\
\partial_{a_{j}}x_{i}^{*}-\partial_{a_{j}}\nu^{*} & =0,\quad\forall j\neq i\label{eq:kkt_ev_a2}\\
\lambda_{i}^{*}\partial_{a_{j}}x_{i}^{*}+x_{i}^{*}\partial_{a_{j}}\lambda_{i}^{*} & =0,\quad\forall j\in[T]\label{eq:kkt_ev_a3}\\
\mu_{i}^{*}\partial_{a_{j}}x_{i}^{*}+(x_{i}^{*}-a_{i})\partial_{a_{j}}\mu_{i}^{*}-\mu_{i}^{*} & =0\label{eq:kkt_ev_a4}\\
\mu_{i}^{*}\partial_{a_{j}}x_{i}^{*}+(x_{i}^{*}-a_{i})\partial_{a_{j}}\mu_{i}^{*} & =0,\quad\forall j\neq i\label{eq:kkt_ev_a5}\\
\sum_{j=1}^{n}\partial_{a_{i}}x_{j}^{*} & =0.\label{eq:kkt_ev_a6}
\end{align}
From (\ref{eq:kkt_ev_a1}) and (\ref{eq:kkt_ev_a2}), we know that
the following holds for all $i$:
\begin{equation}
\partial_{a_{i}}\nu^{*}=\partial_{a_{i}}x_{i}^{*}-\partial_{a_{i}}\lambda_{i}^{*}+\partial_{a_{i}}\mu_{i}^{*}=\partial_{a_{i}}x_{j}^{*},\quad\forall j\neq i.\label{eq:kkt_ev_a7}
\end{equation}
The first equality in~(\ref{eq:kkt_ev_a7}) implies that there exists
a constant $C_{i}$ such that $\partial_{a_{i}}x_{j}^{*}=C_{i}$ for
all $j\neq i$. Then, we can rewrite~(\ref{eq:kkt_ev_a6}) as
\begin{equation}
\partial_{a_{i}}x_{i}^{*}+(T-1)C_{i}=0,\quad\forall i\in[T],\label{eq:kkt_ev_a6-2}
\end{equation}
which implies that 
\begin{align}
 & \norm{\partial_{a_{i}}x^{*}}_{1}=\sum_{j=1}^{T}|\partial_{a_{i}}x_{j}^{*}|=|\partial_{a_{i}}x_{i}^{*}|+\sum_{j\neq i}|\partial_{a_{i}}x_{j}^{*}|\nonumber \\
 & \qquad=2(T-1)|C_{i}|.\label{eq:norm_x_C_i}
\end{align}
Suppose strict complementary slackness holds. Then, for any $i\in[T]$,
only one of the three following cases holds:
\begin{enumerate}
\item $x_{i}^{*}=0$, $\lambda_{i}^{*}>0$, $\mu_{i}^{*}=0$;
\item $x_{i}^{*}=a_{i}$, $\lambda_{i}^{*}=0$, $\mu_{i}^{*}>0$;
\item $0<x_{i}^{*}<a_{i}$, $\lambda_{i}^{*}=0$, $\mu_{i}^{*}=0$.
\end{enumerate}
In the next, we will derive the expression of $\norm{\partial_{a_{i}}x^{*}}_{1}$
for the three cases separately: 
\begin{enumerate}
\item $x_{i}^{*}=0$, $\lambda_{i}^{*}>0$, $\mu_{i}^{*}=0$:

Using~(\ref{eq:kkt_ev_a3}), we obtain $\partial_{a_{j}}x_{i}^{*}=0$
for all $j\in[T]$; in particular, this implies $\partial_{a_{i}}x_{i}^{*}=0$.
Substituting $\partial_{a_{i}}x_{i}^{*}=0$ into~(\ref{eq:kkt_ev_a6-2}),
we obtain $C_{i}=0$, so that $\norm{\partial_{a_{i}}x^{*}}_{1}=0$
according to~(\ref{eq:norm_x_C_i}).

\item $x_{i}^{*}=a_{i}$, $\lambda_{i}^{*}=0$, $\mu_{i}^{*}>0$:

Using~(\ref{eq:kkt_ev_a4}), we obtain $\partial_{a_{i}}x_{i}^{*}=1$;
substitute this into~(\ref{eq:kkt_ev_a6-2}) to obtain $C_{i}=-\frac{1}{T-1}$,
so that we have $\norm{\partial_{a_{i}}x^{*}}_{1}=2$ according to~(\ref{eq:norm_x_C_i}).

\item $0<x_{i}^{*}<a_{i}$, $\lambda_{i}^{*}=0$, $\mu_{i}^{*}=0$:

Using~(\ref{eq:kkt_ev_a3})--(\ref{eq:kkt_ev_a5}), we obtain $\partial_{a_{j}}\lambda_{i}^{*}=\partial_{a_{j}}\mu_{i}^{*}=0$
for all $j\in[T]$; in particular, this implies $\partial_{a_{i}}\lambda_{i}^{*}=\partial_{a_{i}}\mu_{i}^{*}=0$.
Then, using~(\ref{eq:kkt_ev_a7}), we have $\partial_{a_{i}}x_{i}^{*}=C_{i}$;
substitute this into~(\ref{eq:kkt_ev_a6-2}) to obtain $C_{i}=0$,
so that $\norm{\partial_{a_{i}}x^{*}}_{1}=0$ according to~(\ref{eq:norm_x_C_i}).

\end{enumerate}
For case 2 in the above, the fact that $\partial_{a_{i}}x_{j}^{*}=C_{i}=-\frac{1}{T-1}\neq0$
for all $j\neq i$ also implies that $x_{j}^{*}\neq a_{j}$ for all
$j\neq i$. Otherwise, the fact $x_{j}^{*}=a_{j}$ would imply $\mu_{j}^{*}>0$
as a result of strict complementary slackness. According to~(\ref{eq:kkt_ev_a5}),
we have 
\[
\mu_{j}^{*}\partial_{a_{i}}x_{j}^{*}+(x_{j}^{*}-a_{j})\partial_{a_{i}}\mu_{j}^{*}=0,\quad\forall i\neq j,
\]
so that $\mu_{j}^{*}>0$ would imply $\partial_{a_{i}}x_{j}^{*}=0$,
which causes a contradiction. To summarize, there exists at most one
$i\in[T]$ such that $\norm{\partial_{a_{i}}x^{*}}_{1}=2$ (i.e.,
case 2 holds), whereas for other $j\neq i$ we have $\norm{\partial_{a_{j}}x^{*}}_{1}=0$
(i.e., either case 1 or 3 holds). This implies that $\sum_{j=1}^{T}\norm{\partial_{a_{j}}x^{*}}_{1}\leq2$,
which completes the proof.

\subsection{Proof of Lemma~\ref{lem:riemann_integrable-1}\label{sub:appendix-b4}}

Define $\mathcal{I}_{\mu}=\{i\colon\mu_{i}^{*}=0,\: x_{i}^{*}=a_{i}\}$.
If $\mathcal{I}_{\mu}$ is empty, then strict complementary slackness
for the constraint $x\preceq a$ holds. For any $i\in\mathcal{I}_{\mu}$,
we have $\nu^{*}=a_{i}-x_{0,i}$ and $\lambda_{i}^{*}=0$ according
to~(\ref{eq:kkt_ev1}) and~(\ref{eq:kkt_ev3}). For any other $j\notin\mathcal{I}_{\mu}$,
one of the following three cases must hold: (1) $x_{j}^{*}=0$; (2)
$x_{j}^{*}=a_{j}$; (3) $0<x_{j}^{*}<a_{j}$. The last case implies
that $\lambda_{j}^{*}=\mu_{j}^{*}=0$, so that we have $x_{j}^{*}=a_{i}-x_{0,i}+x_{0,j}$,
where $i\in\mathcal{I}_{\lambda}$. Since both $b$ and $x_{0}$ are
fixed, and only one $a_{i}$ among all $i\in[T]$ is allowed to change
due to the constraint imposed by the integration path~$L$, we can
use a similar argument as the one in the proof of Lemma~\ref{lem:riemann_integrable}
to conclude that there are finitely many values of $a$ along $L$
such that the constraint $b=\mathbf{1}^{T}x^{*}$ is satisfied.

\end{document}